\documentclass[12pt]{amsart}
\usepackage[british]{babel}
\usepackage{amsmath,amsfonts,amssymb,amstext,amsbsy,amsopn,amsthm,amsxtra,amscd} 
\usepackage{caption}
\usepackage{color}
\usepackage{enumerate}
\usepackage{enumitem}
\usepackage{geometry}
\usepackage{graphicx}
\usepackage[colorlinks=true]{hyperref}
\usepackage[scr=rsfs]{mathalpha}
\usepackage{mathtools}
\usepackage{pdfpages}
\usepackage{subcaption}
\usepackage{xcolor}
\usepackage{color} 
\usepackage{transparent}
\usepackage{listings}
\usepackage{xcolor}

\numberwithin{equation}{section}

\newtheorem{theorem}{Theorem}[section]
\newtheorem{proposition}[theorem]{Proposition}
\newtheorem{corollary}[theorem]{Corollary}
\newtheorem{lemma}[theorem]{Lemma}
\newtheorem{definition}[theorem]{Definition}
\newtheorem{remark}[theorem]{Remark}
\newtheorem{claim}[theorem]{Claim}

\newtheorem{conjecture}[theorem]{Conjecture}

\DeclareMathOperator{\Ann}{Ann}

\DeclareMathOperator{\nexpl}{NEL}
\DeclareMathOperator{\nexpd}{ONED}
\DeclareMathOperator{\pat}{Pat}
\DeclareMathOperator{\st}{\mathsf{St}}
\DeclareMathOperator{\sst}{\mathsf{Sst}}

\newcommand{\bb}[1]{\mathbb{#1}}
\newcommand{\sob}[1]{\ensuremath{{\mathbin |}\raise-.4ex\hbox{$#1$}}}
\newcommand{\llb}{[\![}
\newcommand{\rrb}{]\!]}

\title[Periodicity from $P_\eta(4,n)\le 4n$]{A Modular Structure Theorem for \\ Minimal Periodic Decompositions and \\ Periodicity of Configurations with $P_\eta(4,n) \le 4n$}

\author[C. F. Colle \& E. Garibaldi]{Cleber Fernando Colle$^1$ \& Eduardo Garibaldi$^2$}

\address{
$^1$ Center for Mathematics, Computing and Cognition, 
Federal University of ABC,
09280-560, Santo Andr\'e, SP, Brazil \\
$^2$ Institute of Mathematics, Statistics and Scientific Computing,
University of Campinas,
13083-859, Campinas, SP, Brazil 
}

\thanks{E.~Garibaldi was partially supported by AMSUD240026, CAPES MathAmSud 88881.985536/2024-01, ANR-22-CE40-3348 THERMOGAMAS, FAPESP 24/04685-2. C.~F.~Colle was partially supported by 
FAPESP 24/15612-6.}
\email{}


\begin{document}
\begin{abstract}
Nivat's conjecture asserts that every two-dimensional configuration
$\eta~\colon~\mathbb{Z}^2~\to~\mathcal{A}$ whose rectangular pattern complexity satisfies
$P_\eta(k,n)~\le~kn $ for some $k,n \in \mathbb{N}$ is periodic.
A theorem of Cyr and Kra \cite{CyrKra16} establishes the
conjecture in the short-rectangle case $P_\eta(k,n)~\le~kn$,  with $ k \le 3 $.
Using the algebraic framework of Kari-Szabados \cite{KariSzabados20}
and recent advances on periodic decompositions and one-sided
nonexpansive directions \cite{Colle23,Colle22}, we extend the
Cyr-Kra result to the case $P_\eta(4,n)\le 4n$: every configuration
satisfying this complexity bound is periodic. The key new ingredient is an intermediate structural theorem of
independent interest: for any non-periodic configuration with low
convex pattern complexity and integer-valued alphabet $ \mathcal{A} $ contained in
$\mathbb{Z}_+$, there exist a configuration $\vartheta$ in the orbit closure
of $\eta$, a $\mathbb{Z}$-minimal periodic decomposition
$\vartheta=\vartheta_1+\cdots+\vartheta_m$, a prime $p\in\mathbb{N}$ with
$\mathcal{A}\subset\llb p \rrb$, and pairs of disjoint half-planes
$U_i,V_i\subset\mathbb{Z}^2$ such that the reductions modulo~$p$ of the
components $\vartheta_i$ are fully periodic on $U_i$ and on $V_i$
simultaneously, for each $1\le i\le m$.
\end{abstract}

\maketitle 

\section{Introduction and statement of results}
\label{sec:intro}

\subsection{Background and motivation}
\label{subsec:background}

Let $\mathcal A$ be an alphabet with at least two elements.  A
\emph{configuration} is an element $\eta\in\mathcal A^{\mathbb Z^2}$, i.e., a
two-dimensional bi-infinite array with values in $\mathcal A$.  For
$k,n\in\mathbb N$, the \emph{rectangular complexity} $P_\eta(k,n)$ counts
the number of distinct $k \times n$ rectangular patterns occurring in
$\eta$. A configuration is \emph{periodic} if there exists
$h\in(\mathbb Z^2)^*$ with $\eta_{g+h}=\eta_g$ for all
$g\in\mathbb Z^2$, and \emph{fully periodic} or \emph{globally periodic} 
when two linearly independent periods exist. From now on, we assume that $\mathcal{A}$ is finite.

More than eight decades ago, Morse and Hedlund \cite{MorseHedlund38,MorseHedlund42} proved that a
one-dimensional sequence is periodic if and only if its complexity
function is eventually bounded by the identity function. Proposed
by Maurice Nivat in his invited talk at ICALP 1997 in Bologna
\cite{Nivat97}, the following conjecture extends the Morse-Hedlund
theorem to two dimensions.

\begin{conjecture}[Nivat, 1997]\label{conj:Nivat}
  Let $\eta\in\mathcal A^{\mathbb Z^2}$.  If $P_\eta(k,n)\le kn$ for some
  $k,n\in\mathbb N$, then $\eta$ is periodic.
\end{conjecture}

\noindent
Conjecture~\ref{conj:Nivat} remains open despite substantial progress
over the past two decades. We briefly recall the milestones most
relevant to the present work.

\medskip\noindent\textit{First partial results.}
Sander and Tijdeman \cite{SanderTijdeman02} proved the conjecture for
$k=2$: if $P_\eta(2,n)\le 2n$, then $\eta$ is periodic. Epifanio, Koskas and Mignosi \cite{EpifanioKosasMignosi03} handled the
case $P_\eta(k,n) \le kn/144$, and Quas and Zamboni \cite{QuasZamboni04} the case $P_\eta(k,n)\le kn/16$.

\medskip\noindent\textit{The Cyr-Kra approach via nonexpansive subdynamics.}
Boyle and Lind \cite{BoyleLind97} introduced the notion of
\emph{expansive subspaces} and proved that any infinite subshift of
$\mathcal A^{\mathbb Z^2}$ admits at least one \emph{nonexpansive line}.
Cyr and Kra \cite{CyrKra15} refined this notion to
\emph{one-sided nonexpansive directions} and used the resulting
geometric structure to prove Nivat's conjecture for the situation
$P_\eta(k,n)~\le~kn/2$.
Their subsequent paper \cite{CyrKra16} establishes the conjecture in
the short-rectangle regime: if $P_\eta(3,n)\le 3n$ for some $n$,
then $\eta$ is periodic. Colle and Garibaldi \cite{ColleGaribaldi20} obtained an alphabetical refinement of the Cyr-Kra bound: for configurations containing all letters of $\mathcal A$, the weaker condition $P_\eta(k,n) \le \frac{1}{2}kn + |\mathcal A| - 1 $, stated more generally for quasi-regular shapes, already suffices to conclude periodicity. The present work takes the Cyr-Kra short-rectangle result as its starting point, extending it from $P_\eta(3,n)$ to $P_\eta(4,n)$.

\medskip\noindent\textit{Algebraic approach: periodic decompositions.}
Using tools from algebraic geometry, Kari and Szabados
\cite{KariSzabados20} proved that any low pattern complexity
configuration $\eta\in\mathcal A^{\mathbb Z^2}$ (with $\mathcal A\subset\mathbb Z$) can be written
as a finite sum of periodic configurations over $\mathbb Z$:
$\eta=\eta_1+\cdots+\eta_m$, with $\eta_i \in$ $\bb{Z}^{\bb{Z}^2}$ periodic. 
As a consequence, they showed that Nivat's conjecture holds if $P_{\eta}(k,n) \leq kn$ for infinitely many pairs $k,n \in \bb{N}$. Szabados \cite{Szabados18} subsequently proved that Nivat's conjecture
holds for configurations admitting such a decomposition with
$m=2$ components. 

\medskip\noindent\textit{Periodic accumulation points and the repetitive case.}
Combining the existence of a non-trivial annihilator \cite{KariSzabados20}
with a compactness argument that eliminates one-sided directions of
determinism, Kari and Moutot \cite{KariMoutot20,KariMoutot21} proved
that the orbit closure of any low convex pattern complexity configuration
contains a periodic configuration.  Applying a technique of Cyr and Kra
\cite{CyrKra15} to uniformly recurrent configurations then yields
Nivat's conjecture in the repetitive (minimal subshift) case.

\medskip\noindent\textit{Periodic decompositions and nonexpansive directions.}
In \cite{Colle23} it was shown that, under natural hypotheses on the
orders of non-periodic configurations in the orbit closure, both
Szabados's conjecture and a symmetry property of the nonexpansive structure hold: 
every nonexpansive line admits both of its orientations as one-sided nonexpansive 
directions, so that each nonexpansive line contributes a pair of oppositely oriented one-sided nonexpansive directions. Together these two conditions yield a reduction of Nivat's 
conjecture to configurations satisfying them.
In \cite{Colle22} a sharper analysis of minimal periodic
decompositions shows that one may always find a decomposition where
each component is defined on a \emph{finite} alphabet contained in
$\llb p\rrb:=\{0,1,\ldots,p-1\}$ for a prime $p$; this
\emph{minimal writing over a finite alphabet} provides the
compactness tools exploited in the present paper.

\subsection{Main results}
\label{subsec:results}

Our first main result is a structural theorem of independent interest:

\begin{theorem}[Modular double-periodicity structure]
\label{thm:modular}
Let $\eta\in\mathcal A^{\mathbb Z^2}$, with $\mathcal A\subset\mathbb Z_+$, be a non-periodic
configuration with low convex pattern complexity.
Then there exist
\begin{enumerate}[label=(\roman*)]
  \item a configuration $\vartheta \in \overline{Orb \, (\eta)}$,
  \item a $\mathbb Z$-minimal periodic decomposition
        $\vartheta=\vartheta_1+\cdots+\vartheta_m$,
  \item a prime $p\in\mathbb N$ with $\mathcal A\subset\llb p \rrb$,
  \item disjoint half-planes $U_i,V_i\subset\mathbb Z^2$, for
        $1\le i\le m$,
\end{enumerate}
such that, for each $1\le i\le m$, both
$\bar\vartheta_i\sob{U_i}$ and $\bar\vartheta_i\sob{V_i}$ are fully
periodic, where $\bar\vartheta_i$ denotes the reduction of $\vartheta_i$
modulo~$p$.
\end{theorem}

\noindent
The proof of Theorem~\ref{thm:modular} (see Section~\ref{sec:modular})
combines two ingredients:
\begin{enumerate}
     \item The existence, established by combining Cases~1 and~2 in the
      proof of \cite[Lemma~4.6]{Colle23}, of an $(\ell,\ell')$-region
       -- a convex unbounded region whose boundary
      consists of two semi-infinite edges, one directed along $\ell$ and
      one directed along $\ell'$, where both $\ell$ and $\ell'$ are
      nonexpansive lines on $\overline{Orb \, (\eta)}$ -- on which a certain accumulation
      point of the orbit of $\eta$ is \emph{fully periodic};
        \item The minimal writing over a finite alphabet from \cite{Colle22},
        which provides the compactness needed to ``symmetrize'' the
        double periodicity within the orbit closure, yielding, for
        each component $\vartheta_i$, a pair of disjoint half-planes
        $U_i, V_i$ on which the reduction $\bar\vartheta_i$ is fully
        periodic.
\end{enumerate}

Our second main result extends the Cyr-Kra theorem \cite{CyrKra16}
from the case $P_\eta(3,n)\le 3n$ to $P_\eta(4,n)\le 4n$:

\begin{theorem}
\label{thm:main}
Let $\eta\in\mathcal A^{\mathbb Z^2}$ be a configuration.  If $P_\eta(k,n)\le kn$ with $ k \le 4 $
for some $n\in\mathbb N$, then $\eta$ is periodic.
\end{theorem}

The proof of Theorem~\ref{thm:main} relies on Theorem~\ref{thm:modular} together with the balanced-set propagation of Cyr and Kra~\cite{CyrKra15} and a maximality argument for the half-planes it produces. We now explain the logical structure of the argument. By Remark~\ref{rem:order-bound}, the complexity bound $ P_\eta(4,n) \le 4n $ forces the order $m$ of any minimal periodic decomposition to satisfy $m \le 4$. The case $m=1$ is trivial and $m=2$ is settled by Theorem~\ref{main_theor_szabados}.
For $ m \in \{3, 4\} $, Theorem~\ref{thm:modular} supplies an accumulation point $\vartheta$ of the orbit of $\eta$ admitting a minimal periodic decomposition $ \vartheta = \vartheta_1 + \cdots + \vartheta_m $ and, for each component $ \vartheta_i $, a pair of disjoint half-planes on which the reduction of $ \vartheta_i $ modulo an appropriate prime is fully periodic. Assuming these half-planes to be maximal in a sense that neither can be enlarged while preserving full periodicity, the contradiction is reached by showing that the geometry of a $\vartheta$-generating set, relative to a suitably chosen unbounded convex region carved out by the half-planes (see Definition~\ref{def:stripfree}), allows the balanced-set argument to propagate full periodicity of some component beyond its assigned half-plane, violating maximality.

The remainder of the paper is organized as follows.
Section~\ref{sec:prelim} collects the preliminary material used
throughout: subshifts and complexity notions are recalled in
Subsection~\ref{subsec:subshifts}, the algebraic formalism of periodic
decompositions and the geometry of their associated convex supports is
developed in Subsection~\ref{subsec:decomp}, and the interplay between
nonexpansive directions and generating sets is reviewed in
Subsection~\ref{subsec:nonexp}.
Section~\ref{sec:modular} is devoted to the proof of
Theorem~\ref{thm:modular}: after establishing
Proposition~\ref{prop_fully_periodic_on strips}, which is the key
intermediate step, we introduce, for conciseness of statement, the
concept of a $\mathbb{Z}_p$-star configuration
(Definition~\ref{Zp-star definition}) and deduce a rephrased version of the
theorem as a consequence (Theorem~\ref{thm:modular_reformulated}).
Section~\ref{sec:main} carries out the argument outlined in the previous paragraph in detail, introducing along the way the precise geometric vocabulary -- strip-free regions (Definition~\ref{def:stripfree}), and compatibility (Definition~\ref{def:compatible}) -- needed to treat the cases $m=3$ and $m=4$ in a unified way, and concludes the proof of Theorem~\ref{thm:main}.

\section{Preliminaries}
\label{sec:prelim}

We briefly recall the main objects and results that will be used in
the proofs; full details can be found in the references cited.

Throughout this paper, whenever a periodic decomposition
$\eta=\eta_1+\cdots+\eta_m$ is considered over a ring $R$ -- either
$\mathbb{Z}$ or $\mathbb{Z}_p$ -- the components $\eta_i$ are
understood to take values in $R$. The same symbol $\eta_i$ is used
in both cases; the ring is always clear from context, and usually no
additional notational distinction between the two situations is made.

\subsection{Subshifts and complexity}\label{subsec:subshifts}

For $\eta\in\mathcal{A}^{\mathbb{Z}^2}$ and a finite, non-empty set $\mathcal S\subset\mathbb Z^2$, the \emph{$\mathcal S$-complexity} $P_\eta(\mathcal S)$ counts the distinct $\mathcal S$-patterns of $\eta$. By an \emph{$\mathcal{S}$-pattern} of $\eta$ we mean the element of $\mathcal{A}^{\mathcal{S}}$ defined by $g\in\mathcal{S}\mapsto\eta_{g+u}\in\mathcal{A}$ for some fixed $u\in\mathbb{Z}^2$ -- it is a subconfiguration of $\eta$ seen through the window $\mathcal{S}$ as it slides over all positions
of the lattice; the collection of all such patterns is denoted $\mathrm{Pat}(\mathcal{S},\eta)$, so that $P_\eta(\mathcal{S})=|\mathrm{Pat}(\mathcal{S},\eta)|$. The configuration $\eta$ has \emph{low pattern complexity} if $P_\eta(\mathcal S)\le|\mathcal S|$ for some finite, non-empty $\mathcal S$, and \emph{low convex pattern complexity} if $\mathcal S$ can be taken convex (in the sense that $\mathcal{S} = \mathrm{Conv}(\mathcal S)\cap\mathbb Z^2$). The \emph{shift action} of $\mathbb Z^2$ on $ \mathcal A^{\mathbb Z^2}$ is given by $ (T^u \eta)_g = \eta_{g+u} $ for $ u, g \in \mathbb Z^2 $; endowing $ \mathcal A^{\mathbb Z^2} $ with the product of discrete topologies makes it a compact metrizable space.
The orbit closure $X_\eta:=\overline{\{T^u\eta:u\in\mathbb Z^2\}}$ is then a closed shift-invariant subset of $\mathcal A^{\mathbb Z^2}$, hence a subshift, and every configuration in $X_\eta$ has low convex pattern complexity whenever $\eta$ does.

For a convex set $ \mathcal S \subset \mathbb Z^2 $ with $ \mathrm{Conv}(\mathcal S) $ of positive area, an \emph{edge} of $ \mathcal S $ is a maximal collinear subset of the boundary of $\mathrm{Conv}(\mathcal S)$ that contains more than one point of $ \mathcal S $; we denote by $E(\mathcal S)$ the collection of all edges of $\mathcal S$, each inheriting a natural orientation from the positively oriented boundary of $ \mathrm{Conv}(\mathcal S)$. If $|E(\mathcal{S})| < \infty$, our convention endows each finite edge $w \in E(\mathcal{S})$ with well-defined successor and predecessor edges. Throughout the text, we commit a slight abuse of notation by also calling edge its maximal subset of points in $ \mathbb Z^2 $. A \emph{vertex} of $ \mathcal S $ is a point of $ \mathcal S $ at which the boundary of $ \mathrm{Conv}(\mathcal S)$ changes direction. 

The global notion of periodicity admits a natural local counterpart that will be needed throughout the paper.

\begin{definition}
Let $\eta \in \mathcal{A}^{\bb{Z}^2}$ be a configuration and let $\mathcal{U} \subset \bb{Z}^2$ be a non-empty set. We say that \(\eta\sob{\mathcal{U}}\) is \emph{periodic} if there exists a non-zero vector $h \in \bb{Z}^2$, called \emph{period of \(\eta\sob{\mathcal{U}}\)}, such that $g+h \in \mathcal{U}$ and $(T^{h}\eta)_{g} = \eta_{g}$ for all $g \in \mathcal{U}$.
If \(\eta\sob{\mathcal{U}}\) has two periods linearly independent (over $\bb{R}^2$), we say that \(\eta\sob{\mathcal{U}}\) is \emph{fully peri\-odic}.	
\end{definition}

When $ \mathcal U = \mathbb Z^2 $ this recovers the usual notion of periodicity of $ \eta $. A central question in the study of low complexity configurations is whether local periodic behavior -- that is, periodicity on proper convex subsets $ \mathcal U $ -- must propagate to the whole of $\mathbb Z^2$. The analogue of Nivat's conjecture for convex shapes, also open, asserts that low convex pattern complexity alone is enough to force global periodicity:

\begin{conjecture}[Nivat's conjecture for convex shapes]
Given a configuration $\eta \in \mathcal{A}^{\mathbb{Z}^2}$, if $ \eta $ has low convex 
pattern complexity, then $\eta$ is periodic.
\end{conjecture}

A central tool in the study of Nivat's conjecture and its variants is the notion of a \emph{generating set}, introduced by Cyr and Kra~\cite{CyrKra15} and further developed in~\cite{FranksKra20}.
A finite convex set $\mathcal S \subset \mathbb Z^2$ is called an \emph{$\eta$-generating
set} if, for all $x \in X_{\eta}$, knowing $x$ on $\mathcal S \setminus \{g\}$ uniquely determines $x_g$
for each vertex $g \in \mathcal S$; intuitively, $ \mathcal S $ is a shape from which the configuration can be reconstructed locally. 

A line is called \emph{nonexpansive} if configurations in $ X_\eta $ cannot be distinguished by what they look like in any arbitrarily thick strip around the line; roughly speaking, it is a direction along which the dynamics fails to separate points. Precisely, a line $\ell\subset\mathbb R^2$ through the origin is a \emph{nonexpansive
line} on $X_\eta$ if, for every $t>0$, there exist distinct
$x,y\in X_\eta$ agreeing on the $t$-neighbourhood of $\ell$; the set
of such lines is $\nexpl(\eta)$. 
If $ \ell \notin \nexpl(\eta) $, we say that $ \ell $ is an \emph{expansive line} on $ X_\eta $.
Given an orientation $\pmb{\ell}$ of $\ell$, the
oriented line is a \emph{one-sided nonexpansive direction} on $X_{\eta}$ if there
exist $x\ne y$ in $X_\eta$ agreeing on the half-plane to the left of
$\pmb{\ell}$, that is, the closed half-plane whose boundary is
$\ell$, positively oriented by $\pmb{\ell}$; the set of such oriented
lines is $\nexpd(\eta)$. 
Whenever $ \pmb{\ell} \notin \nexpd(\eta) $, we say that $ \pmb{\ell} $ is an \emph{one-sided expansive direction} on $ X_\eta $.

From now on, we will consider the following notation: given a line $\ell \subset \bb{R}^2$, we use $\pmb{\ell}$ to denote the line $\ell$ endowed with a given orientation -- for an oriented line $\pmb{\ell} \subset \bb{R}^2$, we use $-\pmb{\ell}$ to denote the oriented line antiparallel to $\pmb{\ell}$ (i.e., with the opposite orientation) that determines the same points of $\pmb{\ell}$ in $\bb{R}^2$.

Given a non-empty, convex set $ \mathcal S \subset \mathbb Z^2 $ and an oriented line $\pmb{\ell}$, we use $\pmb{\ell}_{\mathcal{S}}$ to denote the oriented line $\pmb{\ell}'$ parallel to $\pmb{\ell}$ such that $\mathcal{S}$ is contained in the half-plane to the left of $\pmb{\ell}'$ and $\mathcal{S} \cap \pmb{\ell}' \neq \emptyset$. In particular, if $\pmb{\ell}$ is parallel to some edge of $\mathcal{S}$, then the set $ \mathcal S \cap \pmb{\ell}_{\mathcal S} $ consists of all integer points of $ \mathcal S $ that lie on this edge. Among all non-empty convex subsets of a given finite convex set satisfying $ P_\eta(\mathcal S) \le |\mathcal S| $, Cyr and Kra~\cite{CyrKra15} extract a smallest one: by \cite[Lemma~2.5]{CyrKra15}, any such smallest set is an $\eta$-generating set, which we call an \emph{extremal $\eta$-generating set} to distinguish it from the generating sets associated with periodic decompositions in the next subsection. 
A key feature of extremal $\eta$-generating sets, exploited throughout \cite[Section~4]{CyrKra15}, is the inequality
\begin{equation}\label{eq:balanced}
P_\eta(\mathcal S) - P_\eta(\mathcal S \setminus \pmb{\ell}_{\mathcal S}) < |\mathcal S \cap \pmb{\ell}_{\mathcal S}|,
\end{equation}
which asserts that the edge $ \mathcal S \cap \pmb{\ell}_{\mathcal S} $ contributes strictly fewer new patterns than its number of integer points.

The geometry of $\eta$-generating sets imposes a strong constraint on the nonexpansive structure: any one-sided nonexpansive direction must be parallel or antiparallel to some edge of every $\eta$-generating set. It is this restriction that underlines most the partial results on Nivat's conjecture obtained to date.

\subsection{Periodic decompositions and their geometry}\label{subsec:decomp}

Following \cite{KariSzabados20}, we embed configurations with integer
alphabet in the ring of formal power series $\mathbb Z\llb X^{\pm1}\rrb$ and
write $\eta$ as $\sum_{g\in\mathbb Z^2}\eta_g X^g$.  A Laurent polynomial
$\psi(X)\in \mathbb Z\llb X^{\pm1}\rrb$ \emph{annihilates} $\eta$ if $\psi\eta=0$.
The set of annihilators is denoted $\Ann_{\bb{Z}}(\eta)$. The notion of annihilator can be defined naturally for configurations taking values in other algebraic structures, such as finite fields (see \cite{Colle22} for details).

A configuration $\eta\in\mathcal A^{\mathbb Z^2}$ (with $\mathcal A\subset \bb{Z}$) has a
non-trivial annihilator whenever it has low pattern complexity
\cite{KariSzabados20}.  Moreover, if
$(X^{h_1}-1)\cdots(X^{h_m}-1)\in\Ann_{\bb{Z}}(\eta)$ for vectors
$h_1,\ldots,h_m\in\mathbb Z^2$ in pairwise distinct directions, then $\eta$
decomposes as $\eta=\eta_1+\cdots+\eta_m$ with $\eta_i \in \bb{Z}^{\bb{Z}^2}$ periodic with
period $h_i$ \cite{KariSzabados20}. Such a decomposition is called a
\emph{$\mathbb Z$-periodic decomposition}; it is \emph{$\mathbb Z$-minimal} (of
\emph{order} $m$) if no decomposition into fewer terms exists.  

Given a Laurent polynomial $ \phi(X) = \sum_i a_i X^{u_i} $ with $ a_i \in \mathbb{Z} $ and $u_i \in \mathbb{Z}^2$, its \emph{support} is the finite set $ \mathrm{Supp}(\phi) = \{u_i : a_i \neq 0 \} \subset \mathbb{Z}^2 $. The \emph{reflected convex support} of $ \phi $ is then defined as $ \mathcal{S}_\phi := \mathrm{Conv}(-\mathrm{Supp}(\phi)) \cap \mathbb{Z}^2$, the lattice points of the convex hull of the reflected support. When $ \psi $ annihilates $ \eta $, the set $\mathcal{S}_\psi$ is itself an $\eta$-generating set. The
following lemma, whose proof is given in full generality --  for $ R = \mathbb Z $ or any finite field $\bb{Z}_p$ -- in \cite[Lemma~2.4]{Colle22}, is stated here for completeness.

\begin{lemma}[Szabados \cite{Szabados18}] \label{lem:support_generating}
Let $\eta \in \mathcal{A}^{\mathbb Z^2} $, with $ \mathcal A \subset R $, be a configuration. If $ \psi \in \Ann_R(\eta)$, then $ \mathcal S_\psi $ is an $\eta$-generating set.
\end{lemma}

The geometry of $ \mathcal S_\varphi $, for $\varphi(X) = (X^{h_1}-1) \cdots (X^{h_m}-1) \in \Ann_{\bb{Z}_p}(\eta)$, encodes the periodic structure of $\eta$ in two complementary ways: its edges record the directions of the periods of components of the associated decomposition, and it rigidly constrains the edge structure of every other annihilator of $\eta$.  
The first of these properties is made precise by the following lemma.

\begin{lemma}[Szabados \cite{Szabados18}]\label{lem_periods_annihilator}
Let $\eta \in \mathcal{A}^{\bb{Z}^2}$, with \(\mathcal{A} \subset \bb{Z}_p\), be a configuration and suppose $\varphi(X) = (X^{h_1}-1) \cdots (X^{h_m}-1) \in \Ann_{\bb{Z}_p}(\eta)$, where \(h_1, \ldots, h_m \in \bb{Z}^2\) are vectors in pairwise distinct directions. 
The following conditions hold:
\begin{enumerate}[label=(\roman*)]\setlength{\itemsep}{5pt}
	\item the reflected convex support $\mathcal{S}_{\varphi} $ has exactly $ 2m $ edges -- that is,
    $ |E(\mathcal{S}_{\varphi})| = 2m $ -- and, for each \(1 \leq i \leq m\), one edge is parallel to \(h_i\) and another is parallel to \(-h_i\);
	\item if $\pmb{\ell} \in \nexpd(\eta)$, then there exists \(1 \leq i \leq m\) such that \(\pmb{\ell}\) is either parallel or antiparallel to the vector $h_i$.
\end{enumerate} 
\end{lemma}

\begin{remark}
The original result of Szabados \cite{Szabados18} is stated for annihilators in $ \mathbb Z\llb X^{\pm1}\rrb $. The version above, with annihilators in $ \mathbb Z_p\llb X^{\pm1}\rrb $, follows from two ingredients. First, \cite[Theorem~1.5]{Colle22} guarantees that for an appropriate prime $ p $ with $ \mathcal A \subset \llb p \rrb $ the reduction modulo $ p $ of a $\mathbb{Z}$-minimal periodic decomposition is a $\mathbb{Z}_p$-minimal periodic decomposition of the same order\footnote{This
\emph{minimal writing over a finite alphabet} is crucial for the
compactness arguments in Section~\ref{sec:modular}.}, so that the original annihilator $(X^{h_1}-1)\cdots(X^{h_m}-1)$ descends to $ \Ann{\bb{Z}_p}(\eta) $. Second, the proof that $\mathcal{S}_\varphi$ is an $\eta$-generating set \cite[Lemma~2.4]{Colle22} holds in particular for $\mathbb{Z}_p$, so that the geometric conclusions~(i) and~(ii) carry over verbatim.
\end{remark}

The second complementary property is that, in the minimal situation, the edge geometry of $ \mathcal{S}_\varphi$ is universal: for every other annihilator $\psi$, the reflected convex support $ \mathcal{S}_\psi$ must contain an edge parallel to each edge of $\mathcal{S}_\varphi$.

\begin{proposition}[Kari and Szabados \cite{KariSzabados20}]\label{prop_geom_convexset}
Let $\eta \in \mathcal{A}^{\bb{Z}^2}$, with $\mathcal{A} \subset \bb{Z}_p$, be a non-periodic configuration with $\bb{Z}_p$-minimal periodic decomposition $\eta = \eta_1+\cdots+\eta_m$, where $h_i$ is a period for $\eta_i$, and $\varphi(X) = (X^{h_1}-1) \cdots$ $(X^{h_{m}}-1)$. If $\psi \in \Ann_{\bb{Z}_p}(\eta)$, then, for every edge $w \in E(\mathcal{S}_{\varphi})$, there is an edge $w' \in E(\mathcal{S}_{\psi})$ parallel to $w$.
\end{proposition}

A proof of Proposition~\ref{prop_geom_convexset} as stated above can be found in \cite[Proposition~2.9]{Colle22}. Note that, by item~(ii) of Lemma~\ref{lem_periods_annihilator}, any additional edges that $ \mathcal S_\psi $ may have beyond those dictated by Proposition~\ref{prop_geom_convexset} must be directed along expansive lines of $ X_\eta $: a one-sided nonexpansive direction cannot be parallel or antiparallel to any vector outside $ \{h_1, \ldots, h_m\} $.

\begin{remark}[{\cite[Remark~2.11]{Colle23}}]
\label{rem:order-bound}
If $\eta=\eta_1+\cdots+\eta_m$ is a $\mathbb Z$-minimal periodic decomposition and $P_\eta(n,k)\le nk$ for some $n,k\in\mathbb N$, then $m\le\min\{n,k\}$. In particular, $P_\eta(4,n)\le 4n$ implies $m\le 4$ (for $n\ge 4$). 
Note that the same bound holds for $\mathbb{Z}_p$-minimal periodic decompositions: by \cite[Theorem~1.5]{Colle22}, there exists a prime $ p $ such that the order of a $\mathbb{Z}_p$-minimal decomposition coincides with that of the underlying $\mathbb{Z}$-minimal one, whence the bound $m \le \min\{n,k\}$ follows at once.
\end{remark}

Combining the algebraic decomposition with balanced-set technique of Cyr and Kra~\cite{CyrKra15}, Szabados proved that any low complexity configuration that can be decomposed into a sum of two periodic configurations is necessarily periodic.

\begin{theorem}[Szabados \cite{Szabados18}]\label{main_theor_szabados}
Let $\eta \in \mathcal{A}^{\bb{Z}^2}$, with $\mathcal{A} \subset \bb{Z}_p$, and suppose $\eta = \eta_1+\eta_2$ is a $\bb{Z}_p$-minimal periodic decomposition. Then $\eta$ does not have low convex pattern com\-ple\-xi\-ty.
\end{theorem}

The original version of this result considers $\mathbb{Z}$-decompositions and rectangular complexity.
The adaptation of the arguments of the proof for $\mathbb{Z}_p$-decompositions and convex shapes is briefly provided below for the reader's convenience.

\begin{proof}
Let $h_1, h_2 \in \bb{Z}^2$ denote periods for $\eta_1$ and $\eta_2$, respectively, and consider the set $\mathcal{T} = \{s h_1+t h_2 : s,t \in [0,1]\} \cap \bb{Z}^2$. By the Pigeonhole Principle, there exist integers $0 \leq i < j \leq P_{\eta}(\mathcal{T})$ such that 
$$\eta\sob{\mathcal{T}+i h_1} = \eta\sob{\mathcal{T}+j h_1} = (T^{(j-i)h_1}\eta)\sob{\mathcal{T}+i h_1}.$$ 
Set $\mathcal{T}' := \mathcal{T}+i h_1$ and $S_{\mathcal{T}'} := \bigcup_{t \in \bb{Z}} (\mathcal{T}'+t h_2)$. Since $h_1$ is a period for $\eta_1$, it follows that $\eta_1\sob{\mathcal{T}'} = T^{(j-i)h_1}\eta_1\sob{\mathcal{T}'}$. Combined with $\eta\sob{\mathcal{T}'} = (T^{(j-i)h_1}\eta)\sob{\mathcal{T}'}$ and $\eta = \eta_1+\eta_2$, this gives $\eta_2\sob{\mathcal{T}'} = T^{(j-i)h_1}\eta_2\sob{\mathcal{T}'}$. Now, for $g \in S_{\mathcal{T}'}$, pick $t \in \bb{Z}$ such that $g+t h_2 \in \mathcal{T}'$; then
$$(\eta_2)_g \stackrel{(a)}{=} (\eta_2)_{g+t h_2} \stackrel{(b)}{=} (\eta_2)_{g + t h_2+(j-i) h_1} \stackrel{(a)}{=} (\eta_2)_{ g +(j-i )h_1},$$ 
where both equalities labeled $(a)$ use the fact that $ h_2 $ is a period for $ \eta_2 $, and equality $ (b) $ uses $\eta_2\sob{\mathcal{T}'} = T^{(j-i)h_1}\eta_2\sob{\mathcal{T}'}$ evaluated at $ g + t h_2 \in \mathcal T' $. Hence, $\eta_g = \eta_{g +(j-i) h_1}$ for all $ g \in S_{\mathcal T'} $, which proves that $\eta\sob{S_{\mathcal{T}'}} = (T^{(j-i)h_1}\eta)\sob{S_{\mathcal{T}'}}$.

Suppose, by contradiction, that $\eta$ has low convex pattern complexity. 
Since $ \eta = \eta_1 + \eta_2 $ is a $\mathbb{Z}_p$-minimal decomposition, the annihilator $ (X^{h_1} - 1)(X^{h_2} - 1) $ and Lemma~\ref{lem_periods_annihilator} imply that $\nexpl(\eta)$ consists of at most the two lines $ \ell_1 $ and $ \ell_2 $ containing $ h_1 $ and $ h_2 $, respectively. By~\cite[Proposition~2.7]{Colle22}, each $ \ell_i $ satisfies $ \pm \pmb{\ell}_i \in \nexpd(\eta) $, so that any $\eta$-generating set has edges directed along both $ \ell_1 $ and $ \ell_2 $. In particular, taking $ \ell = \ell_2 $, the local agreement of $ \eta $ and $ T^{(j-i)h_1} \eta $ established on $ S_{\mathcal T'} $ provides the starting data for the balanced-set propagation argument of Cyr and Kra~\cite{CyrKra15}: either a multiple of $ h_1 $ is a global period for $ \eta $, or $ \eta $ is $\ell$-periodic. In either case $ \eta $ is periodic, contradicting the minimality of the decomposition. We refer the reader to \cite[Section~5]{CyrKra15} and \cite[Lemma~4]{Szabados18} for details of this propagation argument.
\end{proof}

\subsection{Nonexpansive directions and generating sets} \label{subsec:nonexp}

The algebraic structure of periodic decompositions interacts with the dynamics of $ X_\eta $ through the geometry of generating sets: as Lemma~\ref{lem_genset_noedge_expas} below makes precise, the edges of any $\eta$-generating set determine which oriented lines can be one-sided nonexpansive directions. By a theorem of Boyle and Lind~\cite{BoyleLind97} every infinite
subshift of $\mathcal A^{\mathbb Z^2}$ possesses at least one nonexpansive line, and
Cyr and Kra~\cite{CyrKra15} made this connection precise: 

\begin{lemma}[Cyr and Kra \cite{CyrKra15}]\label{lem_genset_noedge_expas}
Let $\eta \in \mathcal{A}^{\bb{Z}^2}$ and suppose $\pmb{\ell} \subset \bb{R}^2$ is an oriented line through the ori\-gin and $\mathcal{S} \subset \bb{Z}^2$ is a finite, convex set such that $\mathcal{S} \cap \pmb{\ell}_{\mathcal{S}} = \{g\}$ is $\eta$-generated by $\mathcal{S}$. Then $\pmb{\ell} \not\in \nexpd(\eta)$.
\end{lemma}

In particular, every oriented line in $\nexpd(\eta)$ is parallel to some edge of any $\eta$-generating set whose convex hull has positive area.

By \cite[Proposition~1.3]{Colle23}, for a periodic configuration $\eta$, nonexpansive lines and one-sided nonexpansive directions are perfectly symmetric: $\ell\in\nexpl(\eta)$ if and only if
$-\pmb{\ell},\pmb{\ell} \in \nexpd(\eta)$. This equivalence extends beyond the periodic case, at the cost of a natural hypothesis on the order of configurations in the orbit closure. Indeed, by \cite[Theorem~1.10]{Colle23}, if $\eta$ is a non-periodic low convex pattern complexity configuration for which all non-periodic configurations in $ X_\eta $ have the same order, then two conditions hold simultaneously: Szabados's 
conjecture is satisfied for $\eta$, meaning that every line  
containing a period of some component of a $\mathbb Z$-minimal periodic decomposition of $ \eta $ is a nonexpansive line on $ X_\eta $; and $ \ell \in \nexpl(\eta) $ if and only if $-\pmb{\ell}, \pmb{\ell} \in 
\nexpd(\eta) $, so that each nonexpansive line contributes a pair of oppositely oriented one-sided nonexpansive directions.

The regions introduced below are bounded by lines that interact with the
integer lattice, so it is natural to restrict attention to \emph{rational}
oriented lines -- those containing at least two points of $\mathbb{Z}^2$,
or equivalently whose direction can be represented by a nonzero integer
vector. For such a line $\pmb{\ell}$, its \emph{direction vector}
$v_{\pmb{\ell}}$ is the unique primitive vector in $\mathbb{Z}^2$ pointing
in the direction of the orientation of $\pmb{\ell}$, that is, the shortest
nonzero vector in $\mathbb{Z}^2$ whose direction agrees with that orientation.
The nonexpansive lines that arise in our setting are always rational, as
they contain the period vectors $h_i \in \mathbb{Z}^2$, and so direction
vectors are well defined for all oriented lines considered henceforth.

The propagation of local periodicity across regions bounded by nonexpansive directions is governed by the following notion. Given an oriented line $ \pmb{\ell} \subset \mathbb R^2 $ passing through some point $ q = (q_1, q_2) \in \mathbb Z^2 $, with direction vector $ (a, b) $, the associated \emph{half-plane} is
\[
\mathcal H(\pmb{\ell}) := \{ (g_1,g_2) \in \mathbb Z^2 : - b(g_1 - q_1) + a(g_2 - q_2) \ge 0 \};
\]
the definition is independent of $q$. 
A half-plane $\mathcal{H}(\pmb{\ell})$ is called \emph{rational} if $\pmb{\ell}$ is rational; in that case, its
\emph{slope} is the slope of the direction vector $v_{\pmb{\ell}}$, and two rational half-planes are said to have \emph{distinct slopes} if their bounding lines have distinct underlying directions. For a rational oriented line $\pmb{\ell}$, the \emph{adjacent line} $\pmb{\ell}^{(-)}$ \cite[Definition~2.2]{Colle23} is the unique rational oriented line parallel to $ \pmb{\ell} $, disjoint from $\mathcal H(\pmb{\ell})$, and closest to $\pmb{\ell}$ among all such translates; intuitively, $ \pmb{\ell}^{(-)} $ is the first lattice translate of $ \pmb{\ell} $ in the direction away from $\mathcal H(\pmb{\ell})$.

Given two oriented lines $\pmb{\ell}$ and $\pmb{\ell}'$, in distinct directions and with $\pmb{\ell}$ and $\pmb{\ell}'$ not antiparallel, a convex (positively oriented) unbounded set $\mathcal{K} \subset \bb{Z}^2$ with two semi-infinite edges $w,w' \in E(\mathcal{K})$ is called an \emph{$(\pmb{\ell}, \pmb{\ell}')$-region} if $w$ is parallel to $\pmb{\ell}$, $w'$ is parallel to $\pmb{\ell}'$ and $w$ is a predecessor of $w'$. The directions of the two semi-infinite edges are called the \emph{asymptotic directions} of the region. 

Other useful notion is the concept of envelopment. Let $ \mathcal U \subset \mathbb Z^2$ be a finite convex set with $ \mathrm{Conv}(\mathcal U) $ of positive area. A convex set $ \mathcal T \subset \mathbb Z^2$ (possibly unbounded) is \emph{weakly $E(\mathcal U)$-enveloped} if for every edge $ \varpi \in E(\mathcal T) $ there exists an edge $ \omega \in E(\mathcal U) $ parallel to $ \varpi $ such that $ |\omega \cap \mathcal{U}| \le |\varpi \cap \mathcal T|$; in other words, each edge of $\mathcal T$ has a counterpart in $ E(\mathcal U)$ pointing in the same direction that is no longer, as measured by integer points. When in addition $\mathcal T $ and $ \mathcal U $ have the same number of edges, $ \mathcal T$ is said to be \emph{$E(\mathcal{U})$-enveloped}. 

Given a configuration $\eta$ and a rational oriented line $\pmb{\ell}$, we say that $ \eta $ is \emph{$\ell$-periodic} on a convex set $ \mathcal U \subset \bb{Z}^2 $ if $ \eta \sob{\mathcal U} $ has a period collinear with $v_{\pmb{\ell}}$ (not necessarily parallel). The following lemma describes how $\ell_k$-periodicity (respectively, $\ell_1$-periodicity) on an $(\pmb{\ell}_1, \pmb{\ell}_k)$-region of the form $\mathcal{H}(\pmb{\ell}_1) \cap \cdots \cap \mathcal{H}(\pmb{\ell}_k)$ propagates to the adjacent region obtained by removing the bounding half-plane $\mathcal H(\pmb{\ell}_k)$ (respectively, $\mathcal H(\pmb{\ell}_1)$). It is a key local-to-global tool in the proof of Theorem~\ref{thm:modular}. 

\begin{lemma}\label{lem_ext_periodicity_(l,l')-region}
Let $\eta \in \mathcal{A}^{\bb{Z}^2}$, with \(\mathcal{A} \subset \bb{Z}_p\), be a non-periodic configuration with $\bb{Z}_p$-mi\-ni\-mal periodic decomposition $\eta = \eta_1+\cdots+\eta_m$, periods $ h_i $ for $ \eta_i $, and $\varphi(X) = (X^{h_1}-1) \cdots (X^{h_m}-1)$. Let $\pmb{\ell}_1, \ldots, \pmb{\ell}_k \subset \bb{R}^2$ be oriented lines parallel to the edges of $\mathcal{S}_{\varphi}$, ordered so that the edge parallel to $\pmb{\ell}_{i+1}$ is the successor of the edge parallel to $\pmb{\ell}_{i}$, such that $\mathcal{H}(\pmb{\ell}_1) \cap \cdots \cap \mathcal{H}(\pmb{\ell}_k)$ is an $(\pmb{\ell}_1,\pmb{\ell}_k)$-region weakly $E(\mathcal{S}_{\varphi})$-enveloped and with $k$ edges.
\begin{enumerate}[label=(\roman*)]\setlength{\itemsep}{5pt}
	\item If $\eta$ is $\ell_k$-pe\-ri\-o\-dic on $\mathcal{H}(\pmb{\ell}_1) \cap \cdots \cap \mathcal{H}(\pmb{\ell}_k)$, then $\eta$ and thus $\eta-\eta_k$ are $\ell_k$-pe\-ri\-o\-dic on $\mathcal{H}(\pmb{\ell}_1) \cap \cdots \cap \mathcal{H}(\pmb{\ell}_{k-1})$. 
	
	\item If $\eta$ is $\ell_1$-pe\-ri\-o\-dic on $\mathcal{H}(\pmb{\ell}_1) \cap \cdots \cap \mathcal{H}(\pmb{\ell}_k)$, then $\eta$ and thus $\eta-\eta_1$ are $\ell_1$-pe\-ri\-o\-dic on $\mathcal{H}(\pmb{\ell}_2) \cap \cdots \cap \mathcal{H}(\pmb{\ell}_{k})$.
\end{enumerate}
\end{lemma}

\begin{proof} \

\textit{Item~(ii)} is the symmetric counterpart of
item~(i), and we leave the straightforward adaptation to the reader.

\smallskip

\textit{Item~(i)}. Since the region is weakly $E(\mathcal S_\varphi)$-enveloped, the semi-infinite edge directed along $\pmb{\ell}_k$ is parallel to some edge of $ \mathcal S_\varphi $, which by item~(i) of Lemma~\ref{lem_periods_annihilator} is parallel or antiparallel to $h_j$ for some $ 1\le j \le m$; renaming the components if necessary, we may assume $j=k$, so that $ h_k $ is collinear with $v_{\pmb{\ell}_k}$. Since $\eta_k$ is globally periodic with period $ h_k $, it is periodic with period collinear with $v_{\pmb{\ell}_k}$ on all of $\mathbb Z^2$. Hence the periodicity of $ \eta $ with period collinear with $v_{\pmb{\ell}_k}$ on $\mathcal H(\pmb{\ell}_1)\cap \cdots \cap \mathcal{H}(\pmb{\ell}_{k})$ passes to $ \eta - \eta_k $ on the same region. 

Consider now the annihilator $ \hat\varphi(X) := (X^{h_1} - 1) \cdots \widehat{(X^{h_k} - 1)} \cdots (X^{h_m} - 1) \in \Ann_{\mathbb Z_p}(\eta - \eta_k)$. By item~(i) of Lemma~\ref{lem_periods_annihilator}, the reflected convex support $ \mathcal S_{\hat \varphi}$ has edges parallel or antiparallel to $ h_1, \ldots, \hat{h}_k, \ldots, h_m $; in particular, since the region is weakly $E(\mathcal S_\varphi)$-enveloped and $ h_k $ has been removed, $ \mathcal S_{\hat \varphi} $ has consecutive edges parallel to $\pmb{\ell}_{k-1} $ and $\pmb{\ell}_{k+1}$ but none parallel to $ \pmb{\ell}_k$. This means that $\mathcal S_{\hat \varphi}$ can be positioned so that exactly one of its vertices lies outside the region $\mathcal H(\pmb{\ell}_1)\cap \cdots \cap \mathcal{H}(\pmb{\ell}_{k})$ (see Figure \ref{fig1}). 
\begin{figure}[ht]
	\centering\def\svgwidth{9.8cm}
\begingroup%
  \makeatletter%
  \providecommand\color[2][]{%
    \errmessage{(Inkscape) Color is used for the text in Inkscape, but the package 'color.sty' is not loaded}%
    \renewcommand\color[2][]{}%
  }%
  \providecommand\transparent[1]{%
    \errmessage{(Inkscape) Transparency is used (non-zero) for the text in Inkscape, but the package 'transparent.sty' is not loaded}%
    \renewcommand\transparent[1]{}%
  }%
  \providecommand\rotatebox[2]{#2}%
  \newcommand*\fsize{\dimexpr\f@size pt\relax}%
  \newcommand*\lineheight[1]{\fontsize{\fsize}{#1\fsize}\selectfont}%
  \ifx\svgwidth\undefined%
    \setlength{\unitlength}{249.44879727bp}%
    \ifx\svgscale\undefined%
      \relax%
    \else%
      \setlength{\unitlength}{\unitlength * \real{\svgscale}}%
    \fi%
  \else%
    \setlength{\unitlength}{\svgwidth}%
  \fi%
  \global\let\svgwidth\undefined%
  \global\let\svgscale\undefined%
  \makeatother%
  \begin{picture}(1,0.45454567)%
    \lineheight{1}%
    \setlength\tabcolsep{0pt}%
    \put(0,0){\includegraphics[width=\unitlength,page=1]{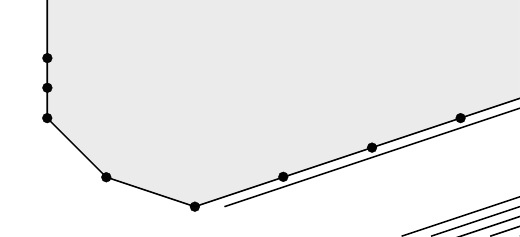}}%
    \put(0.48009858,0.18884447){\color[rgb]{0,0,0}\makebox(0,0)[lt]{\lineheight{1.25}\smash{\begin{tabular}[t]{l}$\mathcal{T}$\end{tabular}}}}%
    \put(0.04359682,0.38943653){\color[rgb]{0,0,0}\makebox(0,0)[lt]{\lineheight{1.25}\smash{\begin{tabular}[t]{l}$\pmb{\ell}_1$\end{tabular}}}}%
    \put(0.01360793,0.24124142){\color[rgb]{0,0,0}\makebox(0,0)[lt]{\lineheight{1.25}\smash{\begin{tabular}[t]{l}$\pmb{\ell}_2$\end{tabular}}}}%
    \put(0.01962119,0.12613535){\color[rgb]{0,0,0}\makebox(0,0)[lt]{\lineheight{1.25}\smash{\begin{tabular}[t]{l}$\pmb{\ell}_3$\end{tabular}}}}%
    \put(0,0){\includegraphics[width=\unitlength,page=2]{fig1.pdf}}%
    \put(0.17310788,0.01817039){\color[rgb]{0,0,0}\makebox(0,0)[lt]{\lineheight{1.25}\smash{\begin{tabular}[t]{l}$\pmb{\ell}_4$\end{tabular}}}}%
    \put(0,0){\includegraphics[width=\unitlength,page=3]{fig1.pdf}}%
  \end{picture}%
\endgroup%

	\caption{A suitable translate $\mathcal{T}$ of $\mathcal{S}_{\hat{\varphi}}$
straddling the boundary of an $(\pmb{\ell}_1,\pmb{\ell}_4)$-region. 
Because $\mathcal{S}_{\hat\varphi}$ has no edge directed along $\ell_4$, 
but does have consecutive edges directed along $ \ell_3$ and 
$ \ell_1 $, it can be positioned so that exactly one vertex lies 
outside the region while all remaining points lie inside, where $\eta-\eta_4$ 
is already known to be $\ell_4$-periodic.}
	\label{fig1}
\end{figure}
Since $ \mathcal S_{\hat \varphi} $ is an $(\eta - \eta_k)$-generating set by Lemma~\ref{lem:support_generating}, the values of $ \eta - \eta_k $ on the points of $ \mathcal S_{\hat \varphi} $ in  
$\mathcal H(\pmb{\ell}_1)\cap \cdots \cap \mathcal{H}(\pmb{\ell}_{k})$ uniquely determine the value at the vertex on the exterior. Starting from the region $\mathcal H(\pmb{\ell}_1)\cap \cdots \cap \mathcal{H}(\pmb{\ell}_{k})$ where $ \eta - \eta_k $ is already known to be periodic with period collinear with $v_{\pmb{\ell}_k}$, we translate  $ \mathcal S_{\hat \varphi} $ successively in the direction of $ v_{\pmb{\ell}_k}$, extending this periodicity of $\eta - \eta_k$ one discrete hal-line at a time across the boundary determined by $ \pmb{\ell}_k$. After countably many such steps, $ \eta - \eta_k $ is periodic with period collinear with $v_{\pmb{\ell}_k}$ on all of $\mathcal H(\pmb{\ell}_1)\cap \cdots \cap \mathcal{H}(\pmb{\ell}_{k-1})$, and hence so is $ \eta $, since $ \eta_k $ is globally periodic with period $ h_k $.
\end{proof}

The following lemma encodes a rigidity phenomenon: if a configuration
is already periodic in some direction, and the configuration seen
asymptotically along that same direction is fully periodic, then full
periodicity propagates back to the entire region.

\begin{lemma}\label{lem:accumulation_fully_periodic}
Let $\xi \in \mathcal{A}^{\mathbb{Z}^2}$ be a configuration and let
$\mathcal{R} \subset \mathbb{Z}^2$ be an $(\pmb{\ell},\pmb{\ell}')$-region.
Let $\pmb{\ell}''$ be a rational oriented line neither parallel nor antiparallel to any semi-infinite edge of $\mathcal{R}$, with $v_{\pmb{\ell}''}$ pointing inward into the region $\mathcal{R}$. Suppose that $\xi$ is $\ell''$-periodic on 
$\mathcal{R}$, and that $\{T^{t v_{\pmb{\ell}''}} \xi : t \in \bb{N}\}$ has a fully periodic accumulation point. Then $\xi$ is fully periodic on $\mathcal{R}$. Consequently, if $\xi$ is $\ell''$-periodic (on the entire lattice), then $\xi$ is fully periodic (on the entire lattice).
\end{lemma}

\begin{proof}
Suppose that $\xi$ is $\ell''$-periodic with period $u''$ on $\mathcal{R} $.
Let $\vartheta$ denote a fully periodic accumulation point of $\{T^{t v_{\pmb{\ell}''}} \xi : t \in \bb{N}\}$
with periods $u''$ along $\ell''$ and $u'$ along $\ell'$.
Let $g \in \mathcal{R}$ be arbitrary. Since $\pmb{\ell}'' \cap \mathrm{int}(\mathrm{Conv}(\mathcal{R}))$ is a
half-line, there exists $t_0 \in \bb{N}$ large enough so that 
\[
\xi\sob{(\{g ,\, g + u' \}+ t_0 u'')}
= \vartheta\sob{\{g,\, g + u'\}}.
\]
Since $\xi$ is $\ell''$-periodic,
\[
\xi_g = \xi_{g + t_0 u''}
\qquad\text{and}\qquad
\xi_{g + u'} = \xi_{g + u' + t_0 u''}.
\]
Since $\vartheta$ is periodic with period $u'$, we have
$\vartheta_{g+u'} = \vartheta_g$, and chaining the equalities above
yields
\[
\xi_{g+u'}
= \xi_{g + u' + t_0 u''}
= \vartheta_{g + u'}
= \vartheta_g
= \xi_{g + t_0 u''}
= \xi_g.
\]
As $g \in \mathcal{R}$ was arbitrary, $\xi$ is periodic with period
$u'$ on $\mathcal{R}$. Together with $\ell''$-periodicity, $\xi$ is
thus fully periodic on $\mathcal{R}$.
\end{proof}

The following proposition, a key ingredient in the proof of Theorem~\ref{thm:main}, extracts periodicity from the balanced inequality~\eqref{eq:balanced}. For its statement we need one further piece of terminology from~\cite{CyrKra15}.
Given a configuration \(\eta \in \mathcal{A}^{\bb{Z}^2}\), a rational oriented line $\pmb{\ell} \subset \bb{R}^2$ through the origin, and a finite convex set $B \subset \bb{Z}^2$, the \emph{semi-strip} of $B$ in the direction of $\pmb{\ell}$ is  
\[
\sst_{B}(\pmb{\ell})
  \;:=\;
  \bigcup_{t\in\mathbb{Z}_+}
  \bigl(B+\,t\,v_{\pmb{\ell}}\bigr).
\]
Let $\mathcal{S} \subset \bb{Z}^2$ be a convex set such that $\mathcal{S} \setminus \pmb{\ell}_{\mathcal{S}} \neq \emptyset$. For each $\gamma \in \pat(\mathcal{S} \setminus \pmb{\ell}_{\mathcal{S}},\eta)$, define 
\[
N_{\mathcal{S}}(\pmb{\ell},\gamma,\eta) = \big|\{\gamma' \in \pat(\mathcal{S},\eta) : \gamma'\sob{\mathcal{S} \setminus \pmb{\ell}_{\mathcal{S}}} = \gamma\}\big|,
\] 
which counts the number of $\mathcal S$-patterns of $\eta$ that extend a given $(\mathcal S \setminus \pmb{\ell}_{\mathcal{S}})$-pattern $ \gamma $; the extension is unique precisely when $N_{\mathcal{S}}(\pmb{\ell},\gamma,\eta) = 1$. A configuration $x \in X_\eta$ is said to be \emph{$(\pmb{\ell},\mathcal{S},\eta,+)$-ambiguous} if every $(\mathcal S \setminus \pmb{\ell}_{\mathcal{S}})$-pattern appearing in the semi-strip $ \sst_{\mathcal S \setminus \pmb{\ell}_{\mathcal S}}(\pmb{\ell}) $ extends to at least two distinct $\mathcal S$-patterns, that is,
\begin{equation}\label{ambiguous+}
	 \ N_{\mathcal S}(\pmb{\ell},(T^{t v_{\pmb{\ell}}} x)\sob{(\mathcal S \setminus \pmb{\ell}_{\mathcal S})}, \eta) > 1, \quad \text{ for all } t \in \mathbb Z_+;
\end{equation}
\emph{$(\pmb{\ell},\mathcal{S},\eta,-)$-ambiguity} is defined analogously with $ \mathbb Z_+ $ replaced by $ \mathbb Z_- $. Whenever $ \pmb{\ell} \in \nexpd(\eta)$ and $\mathcal{S} $ is an $\eta$-gene\-ra\-ting set, ambiguous configurations exist.

\begin{proposition}[Cyr and Kra \cite{CyrKra15}]\label{prop_CyrKra}
Let $\eta \in \mathcal{A}^{\bb{Z}^2}$, let $\pmb{\ell} \in \nexpd(\eta)$,
and let $\mathcal{S} \subset \bb{Z}^2$ be an extremal $\eta$-generating
set. Set $p := |\mathcal{S} \cap \pmb{\ell}_{\mathcal{S}}|-1$, and let
$B \subset \mathcal{S} \setminus \pmb{\ell}_{\mathcal{S}}$ be the largest
convex subset consisting of those points $g \in \mathcal{S} \setminus
\pmb{\ell}_{\mathcal{S}}$ that lie on some translate $\pmb{\ell}'$ of
$\pmb{\ell}$ satisfying $|\mathcal{S} \cap \pmb{\ell}'| \geq p$; that
is, $B$ retains precisely the points of $\mathcal{S} \setminus
\pmb{\ell}_{\mathcal{S}}$ whose row, relative to the direction of
$\pmb{\ell}$, is at least as populated as the leading edge $\mathcal{S} \cap \pmb{\ell}_{\mathcal{S}}$
allows. Suppose $B \neq \emptyset$. Then, for any $(\pmb{\ell},\mathcal{S},
\eta,+)$-ambiguous configuration $x \in X_\eta$:
\begin{enumerate}[label=(\roman*)]\setlength{\itemsep}{5pt}
    \item if $p = 1$, the restriction of $x$ to $\sst_{B}(\pmb{\ell})$
    is $\ell$-periodic;
    \item if $p \geq 2$, the restriction of $x$ to
    $\sst_{B}(\pmb{\ell}) + p\,v_{\pmb{\ell}}$ is $\ell$-periodic.
\end{enumerate}
\end{proposition}

Informally, the proposition asserts that the ambiguity of $x$ along $\pmb{\ell}$ -- the fact that patterns on $\mathcal S \setminus \pmb{\ell}_{\mathcal S}$ fail to determine their extension across the edge $ \pmb{\ell}_{\mathcal S} $ -- forces $ x $ to repeat itself in the direction $ \pmb{\ell} $ on a semi-strip that is as wide as $B$ allows. The complexity hypothesis makes this possible: it bounds the total number of ambiguous extensions across $\pmb{\ell}_{\mathcal S}$ by $|\mathcal S \cap\pmb{\ell}_{\mathcal S }|-1$, which ensures that there are too few distinct $B$-patterns to prevent this repetition. We record the proof here for completeness. 

\begin{proof}
Let $\Gamma = \{(T^{t v_{\pmb{\ell}}}x)\sob{B} : t \in \mathbb{Z}_+\}$. Since $\Gamma \subset \mathrm{Pat}(\mathcal{S} \setminus \pmb{\ell}_{\mathcal{S}},\eta)$, 
\[
p \geq P_{\eta}(\mathcal{S})-P_{\eta}(\mathcal{S} \setminus \pmb{\ell}_{\mathcal{S}}) = \sum_{\gamma \in \pat(\mathcal{S} \setminus \pmb{\ell}_{\mathcal{S}},\eta)} \bigl(N_{\mathcal{S}}(\pmb{\ell},\gamma,\eta)-1\bigr) \geq \sum_{\gamma \in \Gamma} \bigl(N_{\mathcal{S}}(\pmb{\ell},\gamma,\eta)-1\bigr) \geq |\Gamma|.
\] 
Let $\Xi \subset B$ be the set consisting of the initial point of $\mathcal{S} \cap \pmb{\ell}'$ with respect to the orientation of $\pmb{\ell}'$, as $\pmb{\ell}'$ ranges over the oriented lines parallel to $\pmb{\ell}$ that meet $B$. By assumption $B \supset \bigcup_{i=0}^{p} (\Xi + i v_{\pmb{\ell}}) $, and therefore 
\begin{equation}\label{Bound on restrictions}
		p \geq |\Gamma| \geq \big|\big\{ (T^{t v_{\pmb{\ell}}}x)\sob{\bigcup_{i=0}^{p}(\Xi+i v_{\pmb{\ell}})} : t \in \bb{Z}_+\big\}\big|.
	\end{equation}
Setting $\xi_t := (T^{t v_{\pmb{\ell}}}x)\sob{\Xi}$, the bound~\eqref{Bound on restrictions} yields $P_{\xi}(p)\leq p$, where by $P_{\xi}(p)$ we mean the number of distinct words of length $p$ occurring in $\xi$. If $p = 1$ the conclusion is immediate; otherwise, the Morse-Hedlund theorem \cite{MorseHedlund38,MorseHedlund42} implies that $ (\xi_{t+p})_{t \in \bb{Z}_+} $ is periodic, which gives the stated $\ell$-periodicity of $ x $ on the appropriate semi-strip.
\end{proof}

The symmetric statement holds for $(\pmb{\ell}, \mathcal S, \eta, -)$-ambiguous configurations: replacing $\mathbb Z_+$ by $\mathbb Z_-$ in~\eqref{ambiguous+} and reversing the orientation of $\pmb{\ell}$ throughout, the same argument yields the $\ell$-periodicity of $x$ on $\sst_{B}(-\pmb{\ell})$ when $|\mathcal S \cap \pmb{\ell}_{\mathcal S}| = 2 $, and on $\sst_{B}(-\pmb{\ell}) - (|\mathcal S \cap \pmb{\ell}_{\mathcal S}|-1) v_{\pmb{\ell}}$ otherwise.

\section{Proof of the modular double-periodicity theorem}
\label{sec:modular}

This section is devoted to the proof of Theorem~\ref{thm:modular}, the modular double-periodicity structure result stated in the introduction. Roughly speaking, the theorem asserts that, for a non-periodic configuration of  low convex pattern complexity admitting a $\mathbb Z_p$-minimal periodic decomposition $ \vartheta = \vartheta_1 + \cdots + \vartheta_m $, each component $ \vartheta_i $ is fully periodic on two disjoint half-planes. The key intermediate result is the following proposition, which establishes the existence of such half-planes under the hypothesis that $ \vartheta $ is already fully periodic on some $ (\pmb{\ell}, \pmb{\ell}')$-region.

\begin{proposition}\label{prop_fully_periodic_on strips}
Let $\vartheta \in \mathcal{A}^{\bb{Z}^2}$, with \(\mathcal{A} \subset \bb{Z}_p\), be a non-periodic configuration. Suppose $\pmb{\ell}, \pmb{\ell}' \in \nexpd(\vartheta)$ and there exists an $(\pmb{\ell},\pmb{\ell}')$-region $\mathcal{K} \subset \mathbb{Z}^2$ -- that is, a convex unbounded region with asymptotic directions $\pmb{\ell}$ and $\pmb{\ell}'$ -- such that $\vartheta \sob{\mathcal{K}}$ is fully periodic. Then, for any $\bb{Z}_p$-minimal periodic decomposition $\vartheta = \vartheta_1+\cdots+\vartheta_m$, the following conditions hold.   
\begin{enumerate}[label=(\roman*), itemsep=5pt]
	\item There exist half-planes $U_1, \ldots, U_m \subset \bb{Z}^2$, each containing $\mathcal{K} $, such that $\vartheta_i\sob{U_i}$ is fully periodic for each $ 1\le i \le m $. 
	
	\item If, in addition, $\vartheta$ has low convex pattern complexity, then there exist also half-planes $V_1, \ldots, V_m \subset \bb{Z}^2$, with $U_i \cap V_i = \emptyset$, such that $\vartheta_i\sob{V_i}$ is fully periodic for each $ 1\le i \le m $.
\end{enumerate}
\end{proposition}

The proof of condition~(i) proceeds by induction on $m$ propagating full periodicity outward from the region $ \mathcal K $ using the local-to-global tools of Section~\ref{sec:prelim}. The proof of condition~(ii) is more delicate: it requires constructing, for each component, a half-plane on the opposite side of the configuration, and relies on a careful analysis of accumulation points of translates of $\vartheta$ in the direction determined by each nonexpansive line. The key tool is Claim~\ref{claim_zeta_k_fully_periodic} below, which identifies conditions under which a component is forced to be fully periodic at an accumulation point, and then transfers this conclusion back to $\vartheta$ itself.

\begin{proof}
Let \(\vartheta = \vartheta_1+\cdots+\vartheta_m\) be a \(\bb{Z}_p\)-minimal periodic decomposition, let $h_i \in
\mathbb{Z}^2$ be a period for $\vartheta_i$ for $1 \leq i \leq m$, and set  $\varphi(X) := (X^{h_1}-1) \cdots (X^{h_m}-1)$. Let $\pmb{\ell}_1,\ldots,\pmb{\ell}_{2m} \subset \bb{R}^2$ be an enumeration of the oriented lines through the origin parallel to the edges of $\mathcal{S}_{\varphi}$, ordered so that the edge parallel to $\pmb{\ell}_{i+1}$ is the successor of the edge parallel to $\pmb{\ell}_{i}$ in the positively oriented boundary of $\mathcal S_\varphi $, with indices taken modulo $2m$. We may assume that $h_i$ is contained in $\ell_i$ for each $1 \leq i \leq m$. 

Since $\pmb{\ell},\pmb{\ell}' \in \nexpd(\vartheta)$ and $\mathcal{S}_{\varphi}$ is a $\vartheta$-generating set by Lemma~\ref{lem:support_generating}, Lemma \ref{lem_genset_noedge_expas} implies that $\pmb{\ell} = \pmb{\ell}_{\iota}$ and $\pmb{\ell}' = \pmb{\ell}_{\iota'}$ for some $\iota,\iota' \in \{1,2,\ldots,2m\}$. Since $\pmb{\ell}$ and $\pmb{\ell}'$ determine an $(\pmb{\ell},\pmb{\ell}')$-region, one has $\iota' \in \{\iota+1, \ldots,\iota+m-1\}$. The following claim shows that $\iota'$ is in fact uniquely determined by $\iota$.

\begin{claim}\label{claim_iota}
We have $\iota' = \iota+m-1$.
\end{claim}

\begin{proof}[Proof of Claim~\ref{claim_iota}]
Suppose for contradiction $\iota' \neq \iota+m-1$.
Let $\phi(X) = \prod_{i \neq \iota+m-1}(X^{h_i}-1)$. Since $\vartheta\sob{\mathcal{K}}$ is fully periodic, the family $\{T^{t v_{\pmb{\ell}_{\iota+m-1}}}\vartheta\}_{t \in \bb{N}}$, and consequently $\{T^{t v_{\pmb{\ell}_{\iota+m-1}}}\phi(X)\vartheta\}_{t \in \bb{N}}$, has a fully periodic accumulation point. Since \(\phi(X)\vartheta\) is $\ell_{\iota+m-1}$-periodic, Lemma~\ref{lem:accumulation_fully_periodic} implies that \(\phi(X)\vartheta\) is fully periodic. Let $ h \in \bb{Z}^2 $ be a period for \(\phi(X)\vartheta\) such that \(h,h_1, \ldots, h_{m}\) are in pairwise distinct directions, and set $\psi(X) = (X^h-1)\phi(X)$. Then $\psi(X) \in \Ann_{\bb{Z}_p}(\vartheta)$, however, thanks to Lemma~\ref{lem_periods_annihilator}(i), $\mathcal{S}_{\psi}$ has no edge either parallel or antiparallel to $h_{\iota+m-1}$, contradicting Pro\-po\-si\-tion~\ref{prop_geom_convexset}.
\end{proof}

Without loss of generality, assume $\iota = 1$, so that $\iota' = m$. By replacing $\mathcal{K}$ by a subset if necessary, we may further assume that $\mathcal{K}$ is weakly $E(\mathcal{S}_{\varphi})$-enveloped with one edge parallel to $\pmb{\ell}_i$ for each $1 \leq i \leq m$. Lastly, by translating $\vartheta$ if needed, we may suppose that $(0,0) \in (\pmb{\ell}_{m-1})_{\mathcal{K}} \cap (\pmb{\ell}_{m})_{\mathcal{K}}$.

\medskip

\noindent\textit{Step 1: Proof of condition~\textup{(i)} by induction on $m$.}

\smallskip

\textit{Base case $m=2$.} Since $\vartheta\sob{\mathcal{K}}$ is both $\ell_2$-periodic and $\ell_1$-periodic, Lemma \ref{lem_ext_periodicity_(l,l')-region}(i) gives that $\vartheta-\vartheta_2 = \vartheta_1$ is $\ell_2$-periodic on $\mathcal{H}((\pmb{\ell}_1)_{\mathcal{K}})$, and hence fully periodic there; symmetrically, item~(ii) of the same lemma gives that $ \vartheta_2$ is fully periodic on $\mathcal{H}((\pmb{\ell}_2)_{\mathcal{K}})$. Setting $U_i = \mathcal{H}((\pmb{\ell}_i)_{\mathcal{K}})$ for $i = 1,2$ concludes the base case. 

\textit{Inductive step.} Suppose condition (i) holds for every non-periodic configuration of \(\bb{Z}_p\)-order at most $m-1$. Since $\vartheta\sob{\mathcal{K}}$ is $\ell_m$-periodic, Lem\-ma~\ref{lem_ext_periodicity_(l,l')-region}(i) implies that $\vartheta-\vartheta_m = \vartheta_1+\cdots+\vartheta_{m-1}$ is $\ell_m$-periodic on $\mathcal{H}((\pmb{\ell}_1)_{\mathcal{K}}) \cap \cdots \cap  \mathcal{H}((\pmb{\ell}_{m-1})_{\mathcal{K}})$. We claim that every accumulation point of $\{T^{t v_{\pmb{\ell}_m}}(\vartheta_1+\cdots+\vartheta_{m-1}) : t \in \bb{N}\}$ is fully periodic. Indeed, if this failed, the orbit closure of $ \vartheta_1+\cdots+\vartheta_{m-1} $ would contain a configuration $\ell_m$-periodic but not fully periodic, forcing $\ell_m \in \nexpl(\vartheta_1+\cdots+\vartheta_{m-1})$; however item~(ii) of Lemma~\ref{lem_periods_annihilator} applied to the annihilator $(X^{h_1}-1) \cdots (X^{h_{m-1}}-1) \in \Ann_{\bb{Z}_p}(\vartheta_1+\cdots+\vartheta_{m-1})$ would then give a con\-tra\-dic\-tion. Lem\-ma~\ref{lem:accumulation_fully_periodic} then yields that 
$$(\vartheta_1+\cdots+\vartheta_{m-1})\sob{\mathcal{H}((\pmb{\ell}_1)_{\mathcal{K}}) \cap \cdots \cap \mathcal{H}((\pmb{\ell}_{m-1})_{\mathcal{K}})}$$ 
is fully periodic. By induction hypothesis, there exist half-planes $\hat{U}_1, \ldots, \hat{U}_{m-1} $ with $\mathcal{K} \subset \hat{U}_i$ such that $\vartheta_i\sob{\hat{U}_i}$ is fully periodic for each $1 \leq i \leq m-1$. A symmetric argument, applying 
Lem\-ma~\ref{lem_ext_periodicity_(l,l')-region}(ii) to $\vartheta-\vartheta_1 = \vartheta_2+\cdots+\vartheta_m$ yields half-planes $\tilde{U}_2, \ldots, \tilde{U}_m $ with $\mathcal{K} \subset \tilde{U}_i$ such that $\vartheta_i\sob{\tilde{U}_i}$ is fully periodic for $2 \leq i \leq m$. Setting $U_i = \hat{U}_i$ for $1 \leq i \leq m-1$ and $U_m = \tilde{U}_m$ completes the inductive step.

We note for later use that, by the minimality of $m$, the bounding edge of each half-plane $U_i$ is parallel $ h_i$ and hence to $\pmb{\ell}_i$. 

\medskip
\noindent\textit{Step 2: A key claim for condition~\textup{(ii)}.}

\smallskip

The following claim is the engine of the proof of condition (ii). It identifies, for a pair of indices $ j \neq k $, conditions on an accumulation point $ \theta = \theta_1 + \cdots + \theta_m $ of translates of $ \vartheta $ that force $ \theta_k $ to be fully periodic, and then transfers this conclusion back to $ \vartheta_k $ on a half-plane on the opposite side of $\mathcal K$.

\begin{claim}\label{claim_zeta_k_fully_periodic}
Let $j,k \in \{1, \ldots, m\}$, with $j \neq k$, and consider $\tau \in \{-1,1\}$. Suppose that for every accumulation point $\theta$ of $\{T^{t\tau v_{\pmb{\ell}_{j}}}\vartheta : t \in \bb{N}\}$ with $\bb{Z}_p$-periodic decomposition $\theta = \theta_1+\cdots+\theta_m$ -- where $\theta_i$ is an accumulation point of $\{T^{t\tau v_{\pmb{\ell}_{j}}}\vartheta_i : t \in \bb{N}\}$ -- the component $\theta_i$ is fully periodic for every $i \neq j$ and $i \neq k$. Then:
\begin{enumerate}[label=(\alph*), itemsep=5pt]
	\item $\theta_k$ is fully periodic;
	\item there exists a half-plane $V_k $ with bounding edge antiparallel to $\pmb{\ell}_k$ such that $\vartheta_k\sob{V_k}$ is fully periodic.
\end{enumerate}
\end{claim}

\begin{proof}[Proof of Claim~\ref{claim_zeta_k_fully_periodic}] \

\textit{Part~(a).}
If $\theta$ is fully periodic, then $\theta_k = \theta - \sum_{i \neq k} \theta_i$ is $\ell_{j}$-periodic and hence fully periodic. We may therefore assume that $\theta$ is not fully periodic. Since $\theta$ has low convex pattern complexity, Theorem~\ref{main_theor_szabados} implies that $\theta = (\sum_{i \neq k} \theta_i)+\theta_k$ is periodic; consequently, its unique non-expansive line must contain every period. Let $\tilde{h}_{j} \in \ell_{j}$ denote a period for $\sum_{i \neq k} \theta_i$ and set $\psi(X) = (X^{\tilde{h}_j}-1)(X^{h_k}-1)$. Since $\psi(X) \in \Ann_{\bb{Z}_p}(\theta)$, Lemma~\ref{lem_genset_noedge_expas} and item~(ii) of Lem\-ma~\ref{lem_periods_annihilator} imply that the unique non-expansive line on $X_\theta$ is either $\ell_{j}$ or $\ell_k$. Were $\theta$ to be $\ell_k$-periodic, then $\theta_{j} = \theta -\sum_{i \neq j} \theta_i$ would be $\ell_k$-periodic and hence fully periodic; since $\theta_{j}$ is an accumulation point of $\{T^{t\tau v_{\pmb{\ell}_j}}\vartheta_j : t \in \bb{N}\}$, this would force $\vartheta_j$ to be fully periodic, contradicting the minimality of $m$. Therefore $\theta$ is $\ell_{j}$-periodic, whence $\theta_k = \theta - \sum_{i \neq k} \theta_i$ is $\ell_j$-periodic and hence fully periodic. 

\textit{Part~(b).}
For each $\alpha \in \bb{N}$, consider the rectangle $$R_{\alpha} = \{s h_{j}+t h_k : s,t \in [-\alpha,\alpha]\} \cap \bb{Z}^2$$ and the strip $St_{\alpha} = \{ g +t v_{\pmb{\ell}_k} : g \in R_{\alpha}, \ t \in \bb{Z}\}$. 

Let $ (t_n) $ be a sequence in $ \mathbb N $ with $\theta_k = \lim_{n \to +\infty} T^{t_n\tau v_{\pmb{\ell}_{j}}}\vartheta_k$. For each $\alpha$, there exists $t_{n_{\alpha}} $ such that $(T^{t_{n_{\alpha}}\tau v_{\pmb{\ell}_{j}}}\vartheta_k)\sob{R_{\alpha}} = \theta_k\sob{R_{\alpha}}$; the $\ell_k$-periodicity of both $\vartheta_k$ and $\theta_k$ then extends this agreement to the full strip: $(T^{t_{n_{\alpha}}\tau v_{\pmb{\ell}_{j}}}\vartheta_k)\sob{St_{\alpha}} = \theta_k\sob{St_{\alpha}}$. 

Suppose for contradiction that no such half-plane $V_k$ exists. For each $\alpha$, let $L_{\alpha} \subset \bb{Z}^2$ be the largest convex region containing $St_{\alpha}$ such that
\begin{enumerate}[label=(\roman*)]\setlength{\itemsep}{5pt}
	\item $(T^{t_{n_{\alpha}}\tau v_{\pmb{\ell}_{j}}}\vartheta_k)\sob{L_{\alpha}} = \theta_k\sob{L_{\alpha}}$, and
	\item $g \in L_{\alpha}$ if and only if $(g+\ell_k) \cap \bb{Z}^2 \subset L_{\alpha}$.
\end{enumerate}
By the contradiction hypothesis, each $L_{\alpha}$ is either a half-plane with bounding edge parallel to $\pmb{\ell}_k$, or a strip. Let $w_{\alpha} \in E(L_{\alpha})$ be the edge of $L_\alpha$ parallel to $\pmb{\ell}_k$, let $\sigma = \mathrm{sgn}(k-j)$, and let $s_{\alpha} \in \bb{N}$ be such that $s_{\alpha} \sigma t_kv_{\pmb{\ell}_{j}}$ is the lattice point of this form closest to $w_{\alpha}$, where $t_k \in \bb{N}$ is taking so that $t_kv_{\pmb{\ell}_{j}}$ is a period for $\theta_k$ (see Figures~\ref{fig2}{\color{red}(A)} and \ref{fig2}{\color{red}(B)}). 
\begin{figure}[ht]
	\begin{minipage}[t]{0.47\linewidth}
		\begin{center}
			\def\svgwidth{6.8cm}
\begingroup%
  \makeatletter%
  \providecommand\color[2][]{%
    \errmessage{(Inkscape) Color is used for the text in Inkscape, but the package 'color.sty' is not loaded}%
    \renewcommand\color[2][]{}%
  }%
  \providecommand\transparent[1]{%
    \errmessage{(Inkscape) Transparency is used (non-zero) for the text in Inkscape, but the package 'transparent.sty' is not loaded}%
    \renewcommand\transparent[1]{}%
  }%
  \providecommand\rotatebox[2]{#2}%
  \newcommand*\fsize{\dimexpr\f@size pt\relax}%
  \newcommand*\lineheight[1]{\fontsize{\fsize}{#1\fsize}\selectfont}%
  \ifx\svgwidth\undefined%
    \setlength{\unitlength}{214.35772993bp}%
    \ifx\svgscale\undefined%
      \relax%
    \else%
      \setlength{\unitlength}{\unitlength * \real{\svgscale}}%
    \fi%
  \else%
    \setlength{\unitlength}{\svgwidth}%
  \fi%
  \global\let\svgwidth\undefined%
  \global\let\svgscale\undefined%
  \makeatother%
  \begin{picture}(1,0.72731462)%
    \lineheight{1}%
    \setlength\tabcolsep{0pt}%
    \put(0,0){\includegraphics[width=\unitlength,page=1]{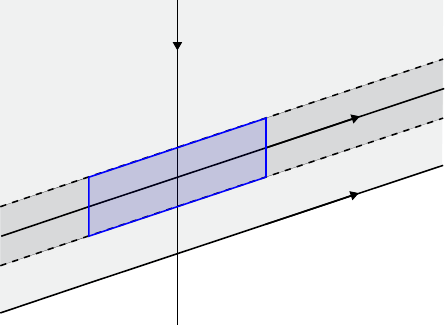}}%
    \put(0.79598122,0.49039089){\color[rgb]{0,0,0}\rotatebox{18.709118}{\makebox(0,0)[lt]{\lineheight{1.25}\smash{\begin{tabular}[t]{l}$\pmb{\ell}_k$\end{tabular}}}}}%
    \put(0,0){\includegraphics[width=\unitlength,page=2]{fig2a.pdf}}%
    \put(0.03509304,0.6557009){\color[rgb]{0,0,0}\makebox(0,0)[lt]{\lineheight{1.25}\smash{\begin{tabular}[t]{l}$L_{\alpha}$\end{tabular}}}}%
    \put(0.79920001,0.25029747){\color[rgb]{0,0,0}\makebox(0,0)[lt]{\lineheight{1.25}\smash{\begin{tabular}[t]{l}$w_{\alpha}$\end{tabular}}}}%
    \put(0.21436755,0.38426152){\color[rgb]{0,0,1}\makebox(0,0)[lt]{\lineheight{1.25}\smash{\begin{tabular}[t]{l}$R_{\alpha}$\end{tabular}}}}%
    \put(0.42520991,0.11041974){\color[rgb]{0,0,0}\makebox(0,0)[lt]{\lineheight{1.25}\smash{\begin{tabular}[t]{l}$s_{\alpha}t_k v_{\pmb{\ell}_j}$\end{tabular}}}}%
    \put(0.4228859,0.6032487){\color[rgb]{0,0,0}\makebox(0,0)[lt]{\lineheight{1.25}\smash{\begin{tabular}[t]{l}$\pmb{\ell}_j$\end{tabular}}}}%
  \end{picture}%
\endgroup%
\\
		\end{center}
		\small{(A) The case $j < k$ ($\sigma = 1$): the maximal region $L_{\alpha}$ lies to one side of its bounding edge $ w_\alpha $, and the black point  $s_{\alpha} t_k v_{\pmb{\ell}_{j}}$ is the translate from the origin by multiples of $ t_k $ in the direction of $\pmb{\ell}_j$ closest to $w_{\alpha}$ from outside $L_\alpha$.}
	\end{minipage} 
	\hfill
	\begin{minipage}[t]{0.47\linewidth}
		\begin{center}
			\def\svgwidth{6.8cm}
\begingroup%
  \makeatletter%
  \providecommand\color[2][]{%
    \errmessage{(Inkscape) Color is used for the text in Inkscape, but the package 'color.sty' is not loaded}%
    \renewcommand\color[2][]{}%
  }%
  \providecommand\transparent[1]{%
    \errmessage{(Inkscape) Transparency is used (non-zero) for the text in Inkscape, but the package 'transparent.sty' is not loaded}%
    \renewcommand\transparent[1]{}%
  }%
  \providecommand\rotatebox[2]{#2}%
  \newcommand*\fsize{\dimexpr\f@size pt\relax}%
  \newcommand*\lineheight[1]{\fontsize{\fsize}{#1\fsize}\selectfont}%
  \ifx\svgwidth\undefined%
    \setlength{\unitlength}{213.84457325bp}%
    \ifx\svgscale\undefined%
      \relax%
    \else%
      \setlength{\unitlength}{\unitlength * \real{\svgscale}}%
    \fi%
  \else%
    \setlength{\unitlength}{\svgwidth}%
  \fi%
  \global\let\svgwidth\undefined%
  \global\let\svgscale\undefined%
  \makeatother%
  \begin{picture}(1,0.72906009)%
    \lineheight{1}%
    \setlength\tabcolsep{0pt}%
    \put(0,0){\includegraphics[width=\unitlength,page=1]{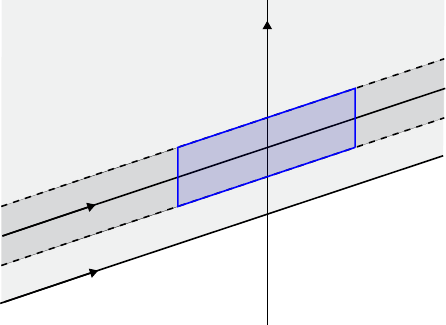}}%
    \put(0.62642589,0.63899529){\color[rgb]{0,0,0}\makebox(0,0)[lt]{\lineheight{1.25}\smash{\begin{tabular}[t]{l}$\pmb{\ell}_j$\end{tabular}}}}%
    \put(0.21234798,0.29962174){\color[rgb]{0,0,0}\rotatebox{18.709118}{\makebox(0,0)[lt]{\lineheight{1.25}\smash{\begin{tabular}[t]{l}$\pmb{\ell}_k$\end{tabular}}}}}%
    \put(0,0){\includegraphics[width=\unitlength,page=2]{fig2b.pdf}}%
    \put(0.03753419,0.65780371){\color[rgb]{0,0,0}\makebox(0,0)[lt]{\lineheight{1.25}\smash{\begin{tabular}[t]{l}$L_{\alpha}$\end{tabular}}}}%
    \put(0.21831678,0.07654625){\color[rgb]{0,0,0}\makebox(0,0)[lt]{\lineheight{1.25}\smash{\begin{tabular}[t]{l}$w_{\alpha}$\end{tabular}}}}%
    \put(0.40809241,0.44973235){\color[rgb]{0,0,1}\makebox(0,0)[lt]{\lineheight{1.25}\smash{\begin{tabular}[t]{l}$R_{\alpha}$\end{tabular}}}}%
    \put(0.62931474,0.22368893){\color[rgb]{0,0,0}\makebox(0,0)[lt]{\lineheight{1.25}\smash{\begin{tabular}[t]{l}$-s_{\alpha}t_k v_{\pmb{\ell}_j}$\end{tabular}}}}%
  \end{picture}%
\endgroup%
\\
		\end{center}
		\small{(B) The case $k<j$ ($\sigma = -1$): the geometry is reversed relative to~(A), and the black point $-s_{\alpha}t_k v_{\pmb{\ell}_{j}}$ is the closest lattice translate to $w_{\alpha}$ in the direction of $-\pmb{\ell}_j$ by multiples of $ t_k $, now lying inside $L_{\alpha}$.}
	\end{minipage}
	\caption{The lattice point $s_\alpha \sigma t_k v_{\pmb{\ell}_j}$ records how far the bounding edge $w_\alpha$ of $L_\alpha$ has drifted from the origin by multiples of $ t_k $ in the direction $\sigma \pmb{\ell}_{j}$, and its sign $ \sigma $ encodes the relative position of $\ell_j$ and $\ell_k$ in the angular ordering of nonexpansive lines.}
	\label{fig2}
\end{figure}
Let $\theta'$ be an accumulation point of $\{T^{(t_{n_{\alpha}}+s_{\alpha} \sigma t_k\tau)\tau v_{\pmb{\ell}_{j}}}\vartheta : \alpha \in \bb{N}\}$; passing to a subsequence, we may assume $\theta_i' = \lim_{l \to +\infty} T^{(t_{n_{\alpha_l}}+s_{\alpha_l} \sigma t_k\tau)\tau v_{\pmb{\ell}_{j}}}\vartheta_i $ exists for each $ i $, giving a $\bb{Z}_p$-periodic decomposition $\theta' = \theta'_1 +\cdots+\theta'_m$. By hypothesis, $\theta'_i$ is fully periodic for every $i \neq j$ and $i \neq k$, so part~(a) implies that $\theta'_k$ is also fully periodic. 

We now show that $ \theta_k = \theta_k' $. Since both configurations are fully periodic, it suffices to verify their agreement on a region large enough to contain a pair of linearly independent periods. 
For each $l \in \bb{N}$, there exists $\alpha_l$ such that $$\theta'_k\sob{R_l} = (T^{(t_{n_{\alpha_l}}+s_{\alpha_l} \sigma t_k\tau)\tau v_{\pmb{\ell}_{j}}}\vartheta_k)\sob{R_l} = (T^{t_{n_{\alpha_l}}\tau v_{\pmb{\ell}_{j}}}\vartheta_k)\sob{(R_l+s_{\alpha_l}\sigma t_k v_{\pmb{\ell}_{j}})}.$$ 
By construction of $ L_{\alpha_l} $, the translated configuration $ T^{t_{n_{\alpha_l}}\tau v_{\pmb{\ell}_{j}}}\vartheta_k $ agrees with $ \theta_k $ on all $ L_{\alpha_l} $, so that $ \theta_k' $ and $T^{s_{\alpha_l}\sigma t_kv_{\pmb{\ell}_j}} \theta_k = \theta_k$ agree on $R_l \cap (L_{\alpha_l} - s_{\alpha_l}\sigma t_kv_{\pmb{\ell}_{j}})$. Since the edge $ w_{\alpha_l} $ lies at bounded distance from $ s_{\alpha_l}\sigma t_kv_{\pmb{\ell}_{j}} $, this intersection grows without bound as $ l \to \infty $. Since $ l $ was arbitrary, we conclude $\theta_k = \theta'_k $. 

Therefore, for any sufficiently large $l \in \mathbb{N}$, the region
$L_{\alpha_l} \cup (St_l + s_{\alpha_l}\sigma t_k v_{\pmb{\ell}_{j}})$ is a
strictly larger convex set than $L_{\alpha_l}$, $\ell_k$-invariant by
construction, on which $T^{t_{n_{\alpha_l}}\tau v_{\pmb{\ell}_{j}}}\vartheta_k$
agrees with $\theta_k$. This contradicts the maximality of $L_{\alpha_l}$
and completes the proof of condition~(b).
\end{proof}

\medskip

\noindent\textit{Step 3: Proof of condition~\textup{(ii)} by descending induction.}

\smallskip

We construct the half-planes $ V_m, V_{m-1}, \ldots, V_1 $ in succession, applying Claim~\ref{claim_zeta_k_fully_periodic} at each step. The strategy is to translate $ \vartheta $
in the direction $ v_{\pmb{\ell}_{p-1}} $ for a descending index $ p $, so that the components 
$ \vartheta_i $ with $ i \neq p-1 $ and $ i \neq p $ are already known to be fully periodic on 
the resulting accumulation points -- either by the half-planes $ U_i $ from Step~1 or by the half-planes
$ V_i $ already constructed -- leaving exactly the pair $ (p-1,p) $ to which the claim applies.

\textit{Construction of $V_m$.} Let $\zeta$ be an accumulation point of $\{T^{t v_{\pmb{\ell}_{m-1}}}\vartheta : t \in \bb{N}\}$; passing to a subsequence, let $\zeta_i = \lim_{n \to +\infty} T^{t_n v_{\pmb{\ell}_{m-1}}}\vartheta_i$ for each $i$, so that $\zeta = \zeta_1 +\cdots+\zeta_m$ is a $\bb{Z}_p$-periodic decomposition. Since $v_{\pmb{\ell}_{m-1}}$ points into the interior of $U_i$ for every $1 \leq i \leq m-2$, the components $\zeta_1,\ldots,\zeta_{m-2}$ are fully periodic. The hypothesis of Claim~\ref{claim_zeta_k_fully_periodic} are therefore satisfied with $ j = m-1 $, $ k = m $, and $ \tau = 1 $, and condition~(b) yields a half-plane $V_m $ with the bounding edge antiparallel to $\pmb{\ell}_m$ such that $\vartheta_m\sob{V_m}$ is fully periodic.

\textit{Descending induction.} Suppose that for some $0 \leq \beta \leq m-3$ the half-planes $ V_m,V_{m-1}, \ldots, $ $ V_{m-\beta}  $ have been constructed with the bounding edge antiparallel to $\pmb{\ell}_{m-i}$, and that  $\vartheta_{m-i}\sob{V_{m-i}}$ is fully periodic for $0 \leq i \leq \beta$. Set $p = m-\beta-1$. Let $\theta$ be an accumulation point of $\{T^{t\pmb{v}_{\pmb{\ell}_{p-1}}}\vartheta : t \in \bb{N}\}$, with $\bb{Z}_p$-periodic decomposition $\theta = \theta_1+\cdots+\theta_m$ obtained by passing to a subsequence. For $p \geq 3$, the vector $ v_{\pmb{\ell}_{p-1}}$ points into  $U_i$ for every $1 \leq i \leq p-2$ and into $V_i$ for every $p+1 \leq i \leq m$; for $p = 2$, it points into $V_i$ for every $p+1 \leq i \leq m$. In both cases, $\theta_i$ is fully periodic for any $i \neq p-1$ and $i \neq p$. Claim~\ref{claim_zeta_k_fully_periodic} with $ j = p-1 $, $ k = p $, and $ \tau = 1$ then yields a half-plane $V_p $ with bounding edge antiparallel to $\pmb{\ell}_p$ such that $\vartheta_p\sob{V_p}$ is fully periodic.

\textit{Construction of $V_1$.} After the descending induction, the half-planes $V_2, \ldots, V_m$ have been defined. To obtain $V_1 $, let $\xi$ be an accumulation point of $\{T^{-t v_{\pmb{\ell}_2}}\vartheta : t \in \bb{N}\}$, with $\bb{Z}_p$-periodic decomposition $\xi = \xi_1+\cdots+\xi_m$ obtained by passing to a subsequence. Since $- v_{\pmb{\ell}_2}$ points into $U_i$ for every $3 \leq i \leq m$, the components $\xi_3, \ldots, \xi_m$ are fully periodic. Applying Claim~\ref{claim_zeta_k_fully_periodic} with $j=2$, $k=1$, and $\tau = -1$ yields a half-plane $V_1$ with bounding edge antiparallel to $\pmb{\ell}_{1}$ such that $\vartheta_{1}\sob{V_1}$ is fully periodic, completing the proof of condition~(ii) and of the proposition.
\end{proof}

We now show how Proposition~\ref{prop_fully_periodic_on strips} implies Theorem~\ref{thm:modular}. To state the conclusion of the theorem in a form amenable to subsequent applications, we introduce the following shorthand.

\begin{definition}\label{Zp-star definition}
Let $p$ be a prime and let $\vartheta \in \mathcal{A}^{\bb{Z}^2}$, with \(\mathcal{A} \subset \llb p \rrb\), be a non-periodic configuration. We say that $\vartheta$ is a \emph{$\bb{Z}_p$-star configuration} if there exist a $\bb{Z}_p$-minimal periodic decomposition $\vartheta = \vartheta_1+\cdots+\vartheta_m$ and rational half-planes $U_1, \ldots, U_m, V_1, \ldots, V_m \subset \bb{Z}^2$, with pairwise distinct slopes and $U_i \cap V_i = \emptyset$ for each $i$, such that $\vartheta_i\sob{U_i}$ and $\vartheta_i\sob{V_i}$ are fully periodic for each $1 \leq i \leq m$.
\end{definition}

Note that the notion of a $\mathbb{Z}_p$-star configuration carries
within itself the full nonexpansive structure of $\vartheta$.
To see this, fix $1 \leq i \leq m$, let $h_i \in \mathbb{Z}^2$ denote
a period for $\vartheta_i$, and let $\ell_i \subset \mathbb{R}^2$ be
the line through the origin containing $h_i$. Clearly, $\ell_i $ is directed 
along the bounding edges of $U_i$ and $ V_i $. Consider any accumulation
point $\theta$ of the sequence $\{T^{t v_{\pmb{\ell}_i}}\vartheta :
t \in \mathbb{Z}\}$: since each component $\vartheta_j$  ($j \neq i$) admits a
multiple of $h_i$ as a period on whichever of $U_j$ or $V_j$ contains
large positive or negative translates of the lattice along $\ell_i$, every component
is eventually invariant under translation by $v_{\pmb{\ell}_i}$ for
$|t|$ sufficiently large, so $\theta$ is $\ell_i$-periodic. On the other
hand, $\theta$ cannot be fully periodic. Indeed, suppose for contradiction
that $\theta$ has a second period in a direction linearly independent
from $h_i$. For $|t|$ sufficiently large, one can find a rectangle
centered at $t v_{\pmb{\ell}_i}$ whose dimensions exceed both $|h_i|$
and the width of the gap between $U_i$ and $V_i$, as well as the norm of a
common period of $\vartheta_j$ ($j \neq i$) in that second direction;
within this rectangle, the translate $T^{tv_{\pmb{\ell}_i}}\vartheta$
inherits the second period, which then propagates to $\vartheta_i$
across a strip containing the gap between $U_i$ and $V_i$, forcing
$\vartheta_i$ to be fully periodic on the entire lattice and contradicting the
minimality of the decomposition. Hence $\theta$ is periodic but not fully
periodic, so as a consequence of the Boyle-Lind theorem~\cite{BoyleLind97} 
$\ell_i$ is its unique nonexpansive line, and from \cite[Proposition~1.3]{Colle23} 
$-\pmb{\ell}_i, \pmb{\ell}_i \in \nexpd(\theta)$. Since
$X_\theta \subset X_\vartheta$, we conclude
$-\pmb{\ell}_i, \pmb{\ell}_i \in \nexpd(\vartheta)$. Repeating for
each $i$ shows that every bounding line contributes both of its
orientations to $\nexpd(\vartheta)$. Conversely, item~(ii) of Lemma~\ref{lem_periods_annihilator} applied
to $\varphi(X)=(X^{h_1}-1)\cdots(X^{h_m}-1)$ shows that every element
of $\nexpd(\vartheta)$ must be either parallel or antiparallel to one of
$h_1,\ldots,h_m$; hence no one-sided nonexpansive direction can arise
from a line other than $\ell_1,\ldots,\ell_m$. Together, the two
inclusions show that $\ell_1,\ldots,\ell_m$ are precisely all the
nonexpansive lines on $X_\vartheta$, each contributing
both of its orientations as one-sided nonexpansive directions.

In language introduced by Definition~\ref{Zp-star definition}, Theorem~\ref{thm:modular} asserts that every non-periodic configuration of low convex pattern complexity has a $\mathbb Z_p$-star accumulation point in its orbit closure:

\begin{theorem}\label{thm:modular_reformulated}
Let $\eta \in \mathcal{A}^{\bb{Z}^2}$, with \(\mathcal{A} \subset \bb{Z}_+\), be a non-periodic configuration with low convex pattern complexity. Then, for some prime $p$ with $\mathcal{A} \subset \llb p \rrb$, there exists a $\bb{Z}_p$-star configuration $\vartheta \in X_\eta$. 
\end{theorem}

The proof proceeds in two steps. The first is to produce, from the orbit closure of $ \eta $, a non-periodic configuration $\vartheta$ that is already fully periodic on some $(\pmb{l},\pmb{l}')$-region $\mathcal K$, with $\pmb{l},\pmb{l}' \in \nexpd(\vartheta)$. This is achieved by revisiting the case analysis of \cite[Lemma~4.6]{Colle23}. The second step applies Proposition~\ref{prop_fully_periodic_on strips} to $\vartheta$: once the existence of such a region is granted, condition~(i) of the proposition yields half-planes $U_i$ on which each component is fully periodic, and condition~(ii) yields the complementary half-planes $V_i$ on the opposite side. The prime $p$ is furnished by \cite[Theorem~1.5]{Colle22}, which guarantees that a $\mathbb Z$-minimal decomposition of $\vartheta$ reduces to a $\mathbb Z_p$-minimal one of the same order.

\begin{proof}
Since $\eta$ has low pattern complexity, Kari and Szabados~\cite{KariSzabados20} guarantee the existence of a non-trivial annihilator. Since $\eta$ is non-periodic, \cite[Theorem~1.9]{Colle23} then yields a line $\ell \in \nexpl(\eta)$ for which $-\pmb{\ell},\pmb{\ell} \in \nexpd(\eta)$. Fix an extremal $\eta$-generating set $\mathcal{S} \subset \bb{Z}^2$, which in particular fulfills inequality~\eqref{eq:balanced}; without loss of generality, we may assume that \(|\mathcal{S} \cap \pmb{\ell}_{\mathcal{S}}| \leq |\mathcal{S} \cap -\pmb{\ell}_{\mathcal{S}}|\). Thus, by \cite[Propositions~2.13 and~2.14]{Colle23}, there exists an $\ell$-periodic configuration $x_{per} \in X_\eta$.
 
The first step is to find a non-periodic configuration $\vartheta \in X_\eta$, oriented lines $\pmb{l}, \pmb{l}' \in \nexpd(\vartheta)$, and a $(\pmb{l},\pmb{l}')$-region $\mathcal{K} $ such that $\vartheta\sob{\mathcal{K}}$ is fully periodic. This is accomplished by the following case analysis, which recapitulates relevant steps of the proof of \cite[Lemma~4.6]{Colle23}. Let \(\eta = \eta_1+\cdots+\eta_m\) be a \(\bb{Z}\)-minimal periodic decomposition, let \(h_i \in \bb{Z}^2\) be a period for \(\eta_i\), and set $\varphi(X) = (X^{h_1}-1) \cdots (X^{h_m}-1)$. 

\medskip\noindent
\textbf{Case 1.} \textit{For some \(E(\mathcal{S}_{\varphi})\)-enveloped set \(B \subset \bb{Z}^2\) there exists $u \in \bb{Z}^2$ such that $(T^{u}\eta)\sob{\sst_B(\pmb{\ell})} = x_{per}\sob{\sst_B(\pmb{\ell})}$.}

\medskip

By an inductive generating-set argument applied to a reduced annihilator -- one from which the factor $ X^{h_k} - 1 $ corresponding to the direction $ \pmb{\ell} $ has been removed -- the agreement of $ T^u \eta $ with $ x_{per} $ on $ \sst_B(\pmb{\ell}) $ propagates $\ell$-periodicity through a nested sequence of regions up to a maximal one, $ \mathcal P $, strictly contained in the full half-plane (otherwise \cite[Proposition~2.14]{Colle23} would force $ T^u \eta $ to be globally periodic, a contradiction; see \cite[Claim~4.7]{Colle23}). The maximality of $ \mathcal P $ then witnesses the ambiguity of $ T^u \eta $ in the direction $ \pmb{\ell}' $ (see \cite[Claim~4.8]{Colle23}), supplying exactly the hypothesis of \cite[Lemma~4.2]{Colle23}. This lemma asserts that a configuration that is already periodic in one asymptotic direction of a region and ambiguous in the other must be globally periodic on some sub-region; applying \cite[Lemma~4.2]{Colle23} here yields an $(-\pmb{\ell},\pmb{\ell}')$-region $ \mathcal K $ on which $T^{u}\eta$ is fully periodic. The minimality of the decomposition then prevents the translates $ \{T^{t v_{\pmb{\ell}'}}(T^u\eta) : t \in \mathbb Z_+ \} $ from having fully periodic accumulation points (see \cite[Claim~4.9]{Colle23}, while any such accumulation point is periodic with period contained in $\ell'$ by \cite[Proposition~2.14]{Colle23}, hence non-fully periodic. The Boyle-Lind theorem and \cite[Proposition~1.3]{Colle23} then yield $-\pmb{\ell}', \pmb{\ell}' \in \nexpd(\eta)$: any such
accumulation point is periodic but not fully periodic, so $\ell'$ is its unique nonexpansive line, and Proposition~1.3 promotes both orientations of $\ell'$ to one-sided nonexpansive directions (see \cite[Claim~4.10]{Colle23}). Since $\nexpd( \cdot )$ is invariant under the shift action, setting $\vartheta = T^{u}\eta$, $\pmb{l} = -\pmb{\ell}$, and $\pmb{l}' = \pmb{\ell}'$ gives the desired data.

\medskip\noindent
\textbf{Case 2.} \textit{For any $E(\mathcal{S}_{\varphi})$-enveloped set $B \subset \bb{Z}^2$ and all $u \in \bb{Z}^2$ with $(T^{u}\eta)\sob{B} = x_{per}\sob{B}$, one has $(T^{u}\eta)\sob{\sst_B(\pmb{\ell})} \neq x_{per}\sob{\sst_B(\pmb{\ell})}$.}

\medskip

Recall that, for a rational oriented line $\pmb{\ell}$, the \emph{adjacent line} $\pmb{\ell}^{(-)}$ \cite[Definition~2.2]{Colle23}
is the rational oriented line parallel to $\pmb{\ell}$, lying just outside $\mathcal{H}(\pmb{\ell})$, closest to it among all such rational lines. The hypothesis of Case~2 is precisely that of \cite[Lemma~3.5]{Colle23}: every shift of $\eta$ that agrees with $ x_{per} $ on some enveloped set $ B $ disagrees with it somewhere in the semi-strip $\sst_B(\pmb{\ell})$. Since $\pmb{\ell} \in \nexpd(\eta)$, there exist configurations in $X_\eta$ that agree on $\mathcal H(\pmb{\ell})$ but differ on $\pmb{\ell}^{(-)}$; in particular, we may assume that $x_{per}\sob{\mathcal{H}(\pmb{\ell}^{(-)})}$ is not fully periodic. \cite[Lemma~3.5]{Colle23} then produces, by a diagonal compactness argument on a nested sequence of enveloped sets exhausting $ \mathcal H(\pmb{\ell}^{(-)}) $, a configuration $ \vartheta \in X_\eta $, an $(\pmb{\ell}, \pmb{\ell}')$-region $\mathcal P$, and an integer $n \in \mathbb{Z}_+$ such that $\vartheta$ agrees with a translate $\hat{x}_{per}$ of $x_{per}$ along $\pmb{\ell}$ on $\mathcal P$ expanded by $n$ layers parallel to $\pmb{\ell}'$, but not on the next layer. By an argument identical to that of Case~1, this maximality witnesses the ambiguity of $\vartheta$ in the direction $\pmb{\ell}'$ (see \cite[Claim~4.12]{Colle23}), and \cite[Lemma~4.2]{Colle23} yields an $(\pmb{\ell},\pmb{\ell}')$-region $\mathcal{K}$, contained in the region of agreement between $\vartheta$ and $\hat{x}_{per}$, on which $\vartheta$ is fully periodic. The conclusion $-\pmb{\ell}', \pmb{\ell}' \in \nexpd(\vartheta) \subseteq \nexpd(\eta)$ then follows by the same reasoning as in Case~1 -- minimality of the decomposition, \cite[Proposition~2.14]{Colle23}, and the Boyle-Lind theorem together with \cite[Proposition~1.3]{Colle23} -- applied now to $\vartheta$ in place of $T^u\eta$ (see \cite[Claims~4.13--4.14]{Colle23}). Finally, $ \vartheta $ is non-periodic: since both $\ell$ and $\ell'$ belong to $\nexpl(\vartheta)$, a periodic but not fully periodic configuration would admit only a single nonexpansive line, a contradiction (see \cite[Claim~4.15]{Colle23}). We therefore set $\pmb{l} = \pmb{\ell}$ and $\pmb{l}' = \pmb{\ell}'$.

\medskip

The two cases differ in their starting point: in Case~1, the driving input is \cite[Claim~4.7]{Colle23}, whose inductive propagation rests on the generating set properties of reduced annihilators -- tools forged by Kari and Szabados~\cite{KariSzabados20}; in Case~2, the corresponding role is played by \cite[Lemma~3.5]{Colle23}, a compactness argument that is conceptually rooted in the balanced-set techniques of Cyr and Kra~\cite{CyrKra15}. Once the respective starting point is in place, the remainder of the argument -- ambiguity, \cite[Lemma~4.2]{Colle23}, and the Boyle-Lind conclusion -- runs along parallel lines, and the two cases converge to the same outcome.

In both cases we obtain a non-periodic configuration $\vartheta \in X_\eta$, directions $\pmb{l},\pmb{l}' \in \nexpd(\vartheta)$, and an $(\pmb{l},\pmb{l}')$-region $\mathcal{K}$ on which $ \vartheta $ is fully periodic. Let $\vartheta = \vartheta_1+\cdots+\vartheta_m$ be any $\bb{Z}$-minimal periodic decomposition. By \cite[Theorem 1.5]{Colle22}, there exists a prime $ p $ with $\mathcal{A} \subset \llb p\rrb$ such that $\overline{\vartheta} = \overline{\vartheta}_1+\cdots+\overline{\vartheta}_m$ is a $\bb{Z}_p$-minimal periodic decomposition of the same order. Proposition~\ref{prop_fully_periodic_on strips}, applied to $ \overline{\vartheta} $ and this decomposition, yields the half-planes $ U_i $ and $ V_i $ required in Definition~\ref{Zp-star definition}, so that $\vartheta$ is a $\bb{Z}_p$-star configuration in $ X_\eta $.
\end{proof}

\section{Proof of the main theorem}
\label{sec:main}

We now prove Theorem~\ref{thm:main}. Suppose that $\eta$ is not periodic. By Theorem~\ref{thm:modular_reformulated}, there exists a $\mathbb Z_p$-star configuration $\vartheta$ in the orbit closure of~$\eta$, admitting a $\mathbb Z_p$-minimal periodic decomposition $\vartheta = \vartheta_1 + \cdots + \vartheta_m $ together with $m$ pairs of complementary rational half-planes $(U_i, V_i)$ with pairwise distinct slopes, such that each component $\vartheta_i$ is fully periodic on both $U_i$ and $V_i$. Since every configuration $\vartheta$ in the orbit closure of $ \eta $ satisfies $P_\vartheta(4,n) \le 4n $, Remark~\ref{rem:order-bound} implies that $ m \le 4$. The case $m=1$ is trivial, and $m=2$ is settled by Theorem~\ref{main_theor_szabados}. The cases $ m = 3$ and $ m = 4 $ are handle in turn below by the same strategy.

In each case the half-planes are assumed to be maximal, meaning that neither $ U_i $ nor $ V_i $ can be enlarged while preserving full periodicity of $\vartheta_i$. The contradiction is reached by showing that one of these half-planes can be extended -- using the geometry of a $ \vartheta$-generating set  and a balanced-set propagation -- so no such maximal star configuration can exist.

To set up the argument, we introduce geometric terminology for collections of complementary half-planes and their intersections.

\begin{definition}\label{def:stripfree}
Let $U = \{U_1, \ldots, U_m\} $ and $V = \{V_1, \ldots, V_m\}$ be collections of rational half-planes in $\bb{Z}^2$ with pairwise distinct slopes and $U_i \cap V_i = \emptyset$ for each~$i$. For $ 1 \le i \le m $, let $\ell_i \subset \bb{R}^2$ denote the line through the origin directed along the bounding edge of $U_i$ (as well as of $V_i$), and let $\st_i(U,V) := \bb{Z}^2 \setminus (U_i \cup V_i)$ be the strip associated with $\ell_i$. We write $(U,V) \in \Pi_m$ to indicate that $ U $ and $ V $ are collections of this form.

Given $(U,V) \in \Pi_m$ and indices $ i \neq j $, an intersection
$$W \in \Big\{ \bigcap_{k = 1}^m Z_k : Z_k \in \{U_k,V_k\}\Big\}$$ 
is called an \emph{$(\pmb{\ell}_i,\pmb{\ell}_j)$-strip-free set} for $(U,V)$ if $W$ is an $(\pmb{\ell}_i,\pmb{\ell}_j)$-region 
in the sense of Subsection~\ref{subsec:nonexp}.

If $\vartheta$ is a $\mathbb Z_p$-star configuration with associated half-planes $(U,V)$, we write $(U,V) \in \Pi_m(\vartheta)$.
\end{definition}

\begin{figure}[ht]
	\centering\def\svgwidth{12.8cm}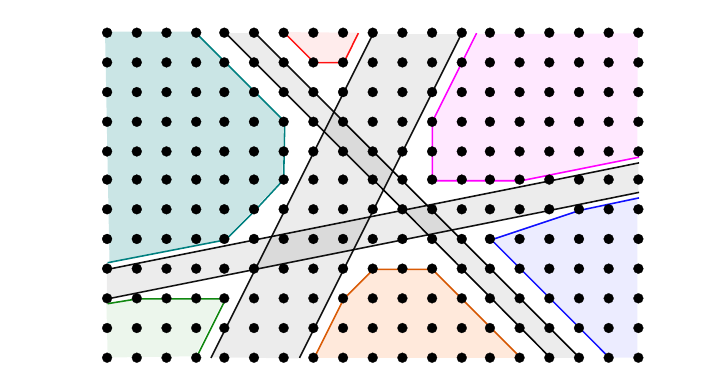
	\caption{An $(\pmb{\ell}_1,\pmb{\ell}_2)$-strip-free set $ A $ for $(U,V) \in \Pi_3$. 
    Replacing all $Z_i$, $ 1 \le i \le 3 $, by the opposite half-planes yields the $(-\pmb{\ell}_1,-\pmb{\ell}_2)$-strip-free $D$.}
	\label{fig3}
\end{figure}

Geometrically, an $(\pmb{\ell}_i,\pmb{\ell}_j)$-strip-free set $ W $ avoids the strips $\st_i(U,V)$ and $\st_j(U,V)$ entirely by lying on a definite side of each of their bounding edges, while the remaining strips at most contribute to determining bounded edges for its boundary. In particular, $W \subset U_k \cup V_k$ for every $1 \leq k \leq m$ (see Figure~\ref{fig3}). 
 
Not every pair of indices $(i,j)$ admits a strip-free set, but at least one adjacent pair always does, as the following remark records. 

\begin{remark}\label{rem:strip-free-existence}
Given $(U,V) \in \Pi_m$ with $m \geq 2$, for each $1 \leq i \leq m$ there exists at least one index $ j \neq i $ such that an $(\pmb{\ell}_i,\pmb{\ell}_j)$-strip-free set for $(U,V)$ exists. Indeed, order the lines $ \ell_1, \dots, \ell_m $ by increasing angle in $[0, \pi)$ and let $ \ell_j $ be the angular predecessor of $\ell_i$ in this ordering, with the convention that the predecessor of the line of smallest angle is the line of largest angle. Endow $\ell_i$ and $\ell_j$ with orientations so that the positively oriented boundary of $ \mathcal H(\pmb{\ell}_i) \cap \mathcal H(\pmb{\ell}_j) $ runs first along the ray in $\ell_i$ and then along the ray in $\ell_j$. Set $Z_i $ to be the element of $ \{U_i, V_i\}$ whose positively oriented boundary is parallel to $\pmb{\ell}_i$, and likewise set $Z_j $ to be the element of $ \{U_j, V_j\}$ whose positively oriented boundary is parallel to $\pmb{\ell}_j$; clearly, $ Z_i \cap Z_j $ is a convex unbounded region. For each remaining index $ k \neq i, j $, choose $ Z_k \in \{U_k, V_k\}$ to be the unique half-plane for which $ Z_i \cap Z_j \cap Z_k $ remains unbounded. The resulting intersection $ W = Z_i \cap Z_j \cap \bigcap_{k \neq i, j} Z_k $ is then an $(\pmb{\ell}_i,\pmb{\ell}_j)$-strip-free set.
Note that $ Z_k $ contributes to determining a bounded edge to the boundary of $ W $ precisely when $ Z_i \cap Z_j $ is not entirely contained in $ Z_k $, as illustrated in Figure~\ref{fig3}.
\end{remark}

The following notion captures the precise geometric condition under which a $\vartheta$-generating set  straddles the boundary of a strip-free region, making it possible to propagate periodicity across that boundary via a balanced-set argument.

\begin{definition}\label{def:compatible}
Given $(U,V) \in \Pi_m$, let $ W = \bigcap_{k=1}^m Z_k $ be an $(\pmb{\ell}_{i},\pmb{\ell}_{j})$-strip-free set for $(U,V)$, and let $\mathcal{T} \subset \bb{Z}^2$ be a finite convex set. Fix $a \in \{i, j\}$, let $b$ denote the other index $(\{b\} = \{i, j\} \setminus \{a\})$, and set $ W^b := \bigcap_{k \ne b} Z_k$, the region obtained from $W$ by dropping the $\ell_b$-constraint; note that $ W^b \supsetneq W $ and that $ W^b $ maintains an unbounded edge parallel to~$\pmb{\ell}_a$.

We say that $W$ is \emph{$\pmb{\ell}_a$-compatible} with $\mathcal T$ if:
\begin{enumerate}[label=(\roman*)]\setlength{\itemsep}{5pt}

    \item $\mathcal{T} \cap \st_{a}(U,V) \neq \emptyset$ and $\sst_{\mathcal{T}}(\sigma^W_a\,\pmb{\ell}_a) \cap \st_{b}(U,V) \neq \emptyset$, where $\sigma^W_a \in \{-1,1\}$ is the sign such that $\sst_{\mathcal{T}}(\sigma^W_a\,\pmb{\ell}_a)$ propagates into the interior of $W$,
	
	\item $\mathcal{T} \setminus (\pmb{\ell}_a)_{\mathcal{T}} \subset W^b$, but $\mathcal{T}  \not\subset W^b$.
\end{enumerate} 
\end{definition}

Condition~(ii) says that all of $\mathcal T$ except its edge $\mathcal{T} \cap (\pmb{\ell}_a)_{\mathcal{T}}$ lies inside the enlarged region $ W^b $, while the edge $\mathcal{T} \cap (\pmb{\ell}_a)_{\mathcal{T}}$ itself protrudes beyond $W^b$; condition~(i) then ensures that this protruding edge does so precisely by reaching into the strip $ \st_a(U, V)$. Together, the two conditions mean that $ \mathcal T $ is positioned so that its edge $\mathcal{T} \cap (\pmb{\ell}_a)_{\mathcal{T}}$ crosses the $\ell_a$-boundary of $ W $, while the remainder of $\mathcal T $ stays safely within $ W^b $.

\begin{remark}\label{rem:edge_inside_halplanes}
Condition~(ii) of Definition~\ref{def:compatible} implies that $\mathcal{T} \cap (\pmb{\ell}_a)_{\mathcal T} \subset Z_k $ for every $ k \neq a $.
Indeed, suppose for contradiction that some point of $\mathcal{T} \cap (\pmb{\ell}_a)_{\mathcal T}$ lies outside $ Z_k $ for some $ k \neq a $. Since the edges of the finite convex set $ \mathcal T $ are arranged in successive angular order along its positively oriented boundary, and since $ \pmb{\ell}_a $ is not directed along $ \ell_k $, this forces the entire edge of $ \mathcal T $ parallel to the bounding edge of $ Z_k $ -- which belongs to $ \mathcal T \setminus (\pmb{\ell}_a)_{\mathcal T} $ -- to lie outside $ Z_k $ as well, contradicting condition~(ii) of Definition~\ref{def:compatible}.
\end{remark}

The contradiction strategy outlined above reduces to a single geometric criterion: given a star configuration $ \vartheta $ with maximal half-planes, it suffices to find a strip-free region whose boundary is straddled by a $\vartheta$-generating set in a controlled way. When such a configuration exists, the balanced-set argument of Cyr and Kra can be initiated, and the maximality assumption is immediately violated. The following theorem makes this precise. 

\begin{theorem}\label{thm:main general theorem}
Let $\vartheta \in \mathcal{A}^{\bb{Z}^2}$ be a $\bb{Z}_p$-star configuration with low convex pattern complexity and maximal half-planes $(U,V) \in \Pi_m(\vartheta)$.
Suppose there exist indices $\imath\neq\jmath$, an $(\pmb{\ell}_{\imath},\pmb{\ell}_{\jmath})$-strip-free set $W $ for $(U,V)$, and a $\vartheta$-generating set $\mathcal T$ such that $ W $ is $\pmb{\ell}_a$-compatible with $ \mathcal{T} $ for some $a \in \{\imath,\jmath\}$. Then:
\begin{enumerate}[label=(\roman*)]\setlength{\itemsep}{5pt}
	\item if $a = \jmath$: $\vartheta$ is $(\pmb{\ell}_a,\mathcal{T},\vartheta,+)$-ambiguous and $\vartheta\sob{\sst_{\mathcal{T} \setminus (\pmb{\ell}_a)_{\mathcal{T}}}(\pmb{\ell}_a)}$ is not $\ell_a$-periodic;
    
    \item if $a = \imath$: $\vartheta$ is $(\pmb{\ell}_a,\mathcal{T},\vartheta,-)$-ambiguous and $\vartheta\sob{\sst_{\mathcal{T} \setminus (\pmb{\ell}_a)_{\mathcal{T}}}(-\pmb{\ell}_a)}$ is not $\ell_a$-periodic.
\end{enumerate} 
\end{theorem}

The proof establishes item~(i) by contradiction: in each case, assuming the contrary allows one to propagate periodicity of a component $ \vartheta_\jmath $ or $ \vartheta_\imath $ beyond the given half-planes, contradicting their maximality. The item~(ii) is entirely symmetric.

\begin{proof}
Let $\vartheta = \vartheta_1+\cdots+\vartheta_m$ denote the $\bb{Z}_p$-minimal periodic decomposition implied in the definition of $\mathbb Z_p$-star configuration (recall Definition~\ref{Zp-star definition}). Since $\vartheta$ has low convex pattern complexity and is non-periodic by definition of~$\mathbb{Z}_p$-star configuration, Theorem~\ref{main_theor_szabados} rules out $m = 2$; hence $m \geq 3$, and in particular every $\vartheta$-generating set has at least $ 2m \geq 6$ edges.

\smallskip\noindent\textit{Item~(i).} 
Set $\pmb{\ell} = \pmb{\ell}_{\imath}$ and $\pmb{\ell}' = \pmb{\ell}_{\jmath}$. Let $h \in \bb{Z}^2$, antiparallel to $ \pmb{\ell} $, be a period vector for $\vartheta_{\imath}$ on the entire lattice, and for each $k \neq \imath$, a period for $ \vartheta_k $ on $  U_k $  (if $h$ points into $U_k$) or on $ V_k $ (if $h$ points into $V_k$). Suppose, for contradiction, that $ \vartheta $ is not $(\pmb{\ell}',\mathcal{T},\vartheta,+)$-ambiguous: then there exists a translate $ \mathcal U \subset \sst_{\mathcal T}(\pmb{\ell}') $ of $\mathcal T$ and a pattern $\gamma = \vartheta \sob{\mathcal U \setminus \pmb{\ell}'_{\mathcal U}}$ such that $ N_{\mathcal{U}}(\pmb{\ell}',\gamma,\vartheta) = 1$. By the choice of the period vector $ \vartheta \sob{(\mathcal U \setminus \pmb{\ell}'_{\mathcal U}) + h} = \vartheta \sob{\mathcal U \setminus \pmb{\ell}'_{\mathcal U}}$, so that the uniqueness of the extension yields
\begin{equation}\label{eq claim main theorem A}
	\vartheta\sob{\mathcal{U}+h} = \vartheta\sob{\mathcal{U}}.
\end{equation} 
Since, by Remark~\ref{rem:edge_inside_halplanes} and periodicity, $ \vartheta_k\sob{\mathcal{U}+h} = \vartheta_k\sob{\mathcal{U}} $ for all $ k \neq \jmath $, equation~\eqref{eq claim main theorem A} yields
$$ (\vartheta_{\jmath})_{g+h} = (\vartheta_{\jmath})_g \qquad \forall \ g \in \pmb{\ell}'_{\mathcal{U}} \cap \mathcal{U}.$$  
To propagate this periodicity along the semi-strip, observe that $ \mathcal U + h + v_{\pmb{\ell}'} \subset W^b $; 
this translate of $ \mathcal T $ is again a $\vartheta$-generating set, so the values of $ \vartheta $ on $ \mathcal U + v_{\pmb{\ell}'} $ excluding the point $ ((\pmb{\ell}'_{\mathcal U} \cap \mathcal{U}) + v_{\pmb{\ell}'}) \setminus (\pmb{\ell}'_{\mathcal U} \cap \mathcal{U}) $ uniquely determine the value of $ \vartheta $ at this vertex. Applying the same reasoning to the translate $\mathcal{U} + v_{\pmb{\ell}'}$ in place of $\mathcal U $, and iterating step by step across $\sst_{\mathcal T}(\pmb{\ell}')$, shows that $ h $ is a period of $\vartheta_{\jmath}$ on an entire semi-strip. This extends the half-plane among $U_{\jmath}$ and $V_{\jmath}$ into whose interior $ h $ points, contradicting the maximality of the half-planes for $\vartheta_{\jmath}$.
 
For the second part, suppose, by contradiction, that $\vartheta\sob{\sst_{\mathcal{T} \setminus \pmb{\ell}'_{\mathcal{T}}}(\pmb{\ell}')}$ is $\ell'$-periodic with period $h' \in \bb{Z}^2$ parallel to $\pmb{\ell}'$, where $ h' $ is assumed to be a period for $ \vartheta_\jmath $ as well. For each $k \neq \jmath $, let $ h' $ also serve as a period for $\vartheta_k $ on $U_k$ or $V_k$ accordingly. Then $\vartheta\sob{W}$ is periodic with period $h'$ and the $\ell'$-periodicity on the semi-strip implies that $\vartheta\sob{W \cup (\pmb{\ell}_W^{(-)} \cap \sst_{\mathcal{T}\setminus \pmb{\ell}'_{\mathcal{T}}}(\pmb{\ell}'))}$ is likewise periodic with period $h'$. By condition~(i) of Definition~\ref{def:compatible}, the semi-strip $\sst_{\mathcal{T}}(\pmb{\ell}')$ meets $\pmb{\ell}_W^{(-)}$, so the intersection $\pmb{\ell}_W^{(-)} \cap
\sst_{\mathcal{T}\setminus(\pmb{\ell}')_{\mathcal{T}}}(\pmb{\ell}')$ is non-empty. Actually, since $ \mathcal T $ has at least $ 2m \ge 6 $ edges, a standard congruence argument for convex lattice sets gives
\begin{equation}\label{eq:cardinality_adjacent}
  \bigl|\pmb{\ell}_W^{(-)} \cap
  \sst_{\mathcal{T}\setminus(\pmb{\ell}')_{\mathcal{T}}}(\pmb{\ell}')
  \bigr|
  \;\geq\;
  \bigl| \pmb{\ell}_{\mathcal{T}} \cap \mathcal{T}\bigr| - 1.
\end{equation}
Let $\mathcal{U}$ be the translate of $\mathcal{T}$ such that the initial points of $(\mathcal{U} \cap \pmb{\ell}_{\mathcal{U}})+v_{\pmb{\ell}}$ and $\pmb{\ell}_W^{(-)} \cap \sst_{\mathcal{T}\setminus(\pmb{\ell}')_{\mathcal{T}}}(\pmb{\ell}')$ coincide. By~\eqref{eq:cardinality_adjacent}, every point of $ \mathcal{U} $ with the exception of  $(\pmb{\ell}_{\mathcal{U}} \cap \mathcal{U}) \setminus ((\pmb{\ell}_{\mathcal{U}} \cap \mathcal{U}) + v_{\pmb{\ell}})$ lies in the region where $\vartheta$ is already known to have $h'$ as a period. Since $\mathcal{U}$ is a $\vartheta$-generating set, the generating property allows to extend the period $h'$ to the vertex $(\pmb{\ell}_{\mathcal{U}} \cap \mathcal{U}) \setminus ((\pmb{\ell}_{\mathcal{U}} \cap \mathcal{U}) + v_{\pmb{\ell}})$ itself: that is, $\vartheta_{g+h'} = \vartheta_g$ for $ \{g\} =  (\pmb{\ell}_{\mathcal{U}} \cap \mathcal{U}) \setminus ((\pmb{\ell}_{\mathcal{U}} \cap \mathcal{U}) + v_{\pmb{\ell}})$. Iterating this extension step by step -- translating $\mathcal{U}$ successively by $-v_{\pmb{\ell}}$ -- propagates the period $ h' $ of $\vartheta$ along an entire semi-ray contained in $\pmb{\ell}_W^{(-)} \cap \mathbb{Z}^2$; this process is illustrated in Figures~\ref{fig4}{\color{red}(A)} and~\ref{fig4}{\color{red}(B)}. Note that this semi-ray is in particular contained in $ W^b $. Therefore, since $h'$ is a period for $\vartheta_k$ ($k \neq \imath$) on the relevant half-planes, it follows that $h'$ is also a period of $\vartheta_\imath$ along this semi-ray, thereby enlarging whichever of $U_\imath$ and $V_\imath$ into whose interior $h'$ points and contradicting their maximality.

\begin{figure}[ht]
	\begin{minipage}[t]{0.47\linewidth}
		\begin{center}
			\def\svgwidth{7.8cm}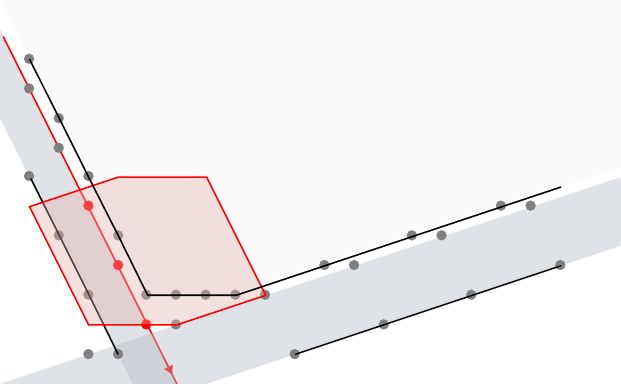\\
		\end{center}
		\small{(A) The red points form the set $\pmb{\ell}_W^{(-)} \cap \sst_\mathcal{T}(\pmb{\ell}_{\jmath})$;  the edge $w$ is parallel to $\pmb{\ell}$ and the edge $w'$ is parallel to $\pmb{\ell}'$.}
	\end{minipage} 
	\hfill
	\begin{minipage}[t]{0.47\linewidth}
		\begin{center}
			\def\svgwidth{7.8cm}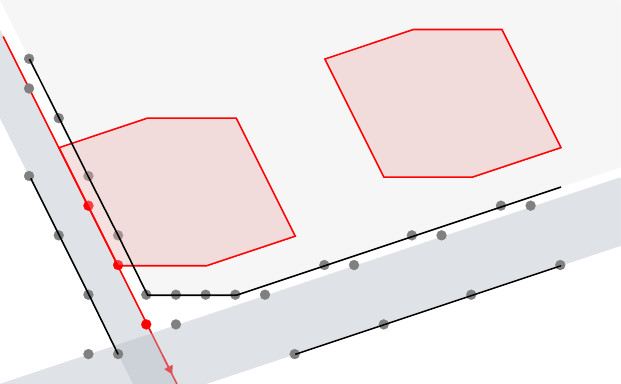\\
		\end{center}
		\small{(B) The set $\mathcal{U}$ is a translation of $\mathcal{T}$; since $\mathcal{U}$ (and any of its translations) is a $\vartheta$-generating set, the values of $\vartheta$ at the blue vertices are uniquely determined by its values at the remaining points of $\mathcal U$ and $\mathcal{U}+h'$, respectively.}
	\end{minipage}
	\caption{Proof of item~(ii): propagation of the period $h'$ of $\vartheta$ along the adjacent line $\pmb{\ell}_W^{(-)}$ via
successive applications of the generating set argument.}
	\label{fig4}
\end{figure}

\smallskip\noindent\textit{Item~(ii)} is the symmetric
counterparts of~(i): replace $\pmb{\ell}_\jmath$ by
$\pmb{\ell}_\imath$ and reverse the relevant orientations throughout.
The argument is identical and we omit the details.
\end{proof}

The theorem above provides the key mechanism for deriving periodicity from
a geometric compatibility condition. We now show how it specialises, together
with the order bound of Remark~\ref{rem:order-bound}, to settle Nivat's
conjecture in the case $P_\eta(4,n) \leq 4n$.

\begin{corollary}
Given $\eta \in \mathcal{A}^{\bb{Z}^2}$, with $\mathcal{A} \subset \bb{Z}_+$,
suppose there exists $n \in \bb{N}$ such that $P_{\eta}(4,n) \leq 4n$.
Then $\eta$ is periodic.
\end{corollary}

\begin{proof}
The argument proceeds by contradiction. Assuming $\eta$ non-periodic,
Theorem~\ref{thm:modular_reformulated} supplies a $\bb{Z}_p$-star
configuration $\vartheta \in X_\eta$ whose order $m$ is either $3$ or $4$,
by Theorem~\ref{main_theor_szabados} and Remark~\ref{rem:order-bound}.
In each case the strategy is the same: we locate a strip-free region $W$
and a translate $\mathcal{T}$ of an extremal $\vartheta$-generating set
$\mathcal{S} \subset R_{4,n} :=  \llb 4 \rrb \times \llb n \rrb $ such that $W$ is $\pmb{\ell}_a$-compatible
with $\mathcal{T}$ for some index $a$. Theorem~\ref{thm:main general theorem}
then gives ambiguity, and Proposition~\ref{prop_CyrKra} converts that
ambiguity into periodicity on a semi-strip, directly contradicting
Theorem~\ref{thm:main general theorem}. The geometry of $\mathcal{S}
\subset R_{4,n}$ -- in particular, the constraints it imposes on edge
cardinalities -- drives the case analysis.

Without further delay, suppose, by contradiction, that $\eta$ is non-periodic. By
Theorem~\ref{thm:modular_reformulated}, for some prime $p$ with
$\mathcal{A} \subset \llb p \rrb$, there exists a $\bb{Z}_p$-star
configuration $\vartheta \in X_{\eta}$ with $\bb{Z}_p$-minimal periodic
decomposition $\vartheta = \vartheta_1+\cdots+\vartheta_m$.
By the very nature of a $\bb{Z}_p$-star configuration (see the comments 
that follow Definition~\ref{Zp-star definition}), a line $\ell$ belongs to
$\nexpl(\vartheta)$ if and only if it contains a period for some
$\vartheta_i$, and $\ell \in \nexpl(\vartheta)$ if and only if
$-\pmb{\ell}, \pmb{\ell} \in \nexpd(\vartheta)$.
By Theorem~\ref{main_theor_szabados} and Remark~\ref{rem:order-bound},
$m \in \{3,4\}$.

Let $(U,V) \in \Pi_m(\vartheta)$ be maximal half-planes and, for each
$1 \le i \le m$, let $\ell_i \subset \bb{R}^2$ denote the line through
the origin directed along the bounding edge of $U_i$ (and of $V_i$).
Let $\mathcal{S} \subset R_{4,n} $ be an extremal $\vartheta$-generating set.

\medbreak

\textbf{Case 1: $m = 3$.}

\medbreak

Any $\vartheta$-generating set contained in $R_{4,n}$ has at most $8$
edges, of which at least $6$ come in antiparallel pairs
(Lemma~\ref{lem_genset_noedge_expas}). 

\medskip

\textit{Case 1(a): the $y$-axis is expansive on $X_{\vartheta}$.}

\medskip

We orient $\pmb{\ell}_1$ so that the bounding edge of $U_1$ is parallel to $\pmb{\ell}_1$. After renaming if necessary, we index the underlying lines $\ell_2$ and $\ell_3$ in the counterclockwise direction starting from $\pmb{\ell}_1$.  The orientations of $\pmb{\ell}_2$ and $\pmb{\ell}_3$ are then chosen so that the following
three strip-free regions are well defined:
\begin{equation}\label{eq:strip-free-triple-m3}
\left.
\begin{aligned}
&W &&\text{an } (-\pmb{\ell}_3,\pmb{\ell}_2)\text{-strip-free region for }(U,V),\\
&W' &&\text{an } (\pmb{\ell}_1,\pmb{\ell}_3)\text{-strip-free region for }(U,V),\\
&W'' &&\text{an } (\pmb{\ell}_2,-\pmb{\ell}_1)\text{-strip-free region for }(U,V).
\end{aligned}
\right.
\end{equation}

Note that no non-vertical edge $w \in E(\mathcal{S})$
can satisfy $|w \cap \mathcal{S}| \geq 3$, with the sole possible exception
of a vertical edge. We therefore ask: does $W'$ admit a translate $\mathcal{T}$ of $\mathcal{S}$ 
with which it is $\pmb{\ell}_3$-compatible? If so, Proposition~\ref{prop_CyrKra}(i) 
yields $\ell_3$-periodicity of $\vartheta$ on 
$\sst_{\mathcal{T} \setminus (\pmb{\ell}_3)_{\mathcal{T}}}(\pmb{\ell}_3)$, 
contradicting Theorem~\ref{thm:main general theorem}.

We may therefore assume that $W'$ is \emph{not} $\pmb{\ell}_3$-compatible
with any translate of $\mathcal{S}$. To exploit this, let $u' \in
(-\pmb{\ell}_3)_{\st_3(U,V)} \cap \bb{Z}^2$ be the unique point with
$u' \notin \st_2(U,V)$ and $u' - v_{\pmb{\ell}_3} \in
\mathcal{H}((-\pmb{\ell}_2)_{\st_2(U,V)})$, and define the parallelogram
\[
\mathcal{Q} := \mathrm{Conv}\!\left\{
u',\;
u'-v_{\pmb{\ell}_3},\;
u'-v_{\pmb{\ell}_2},\;
u'-v_{\pmb{\ell}_3}-v_{\pmb{\ell}_2}
\right\} \cap \bb{Z}^2.
\]
Let $u \in (-\pmb{\ell}_2)_{\st_2(U,V)} \cap \mathcal{Q}$ be such
that $u \notin \st_3(U,V)$, and let $\mathcal{K}$ be the translate
of $\mathcal{S}$ whose final point of the edge parallel to $\pmb{\ell}_2$ coincides with $u$
(see Figure~\ref{fig5}{\color{red}(A)}). Since $W'$ is not $\pmb{\ell}_3$-compatible with any translate of $\mathcal{S}$, necessarily $\mathcal{K}+u'-u \subset U_1$.
Let $ - w_1, w_1 \in E(\mathcal{K}+u'-u)$ denote an antiparallel pair of edges of $\mathcal{K}+u'-u$,
with $ w_1 $ parallel to $\pmb{\ell}_1$. Note that
there exist three vertical lines from the initial vertex of $-w_1$
to the final vertex of $w_1$, inclusive
(see Figure~\ref{fig5}{\color{red}(B)}).

\begin{figure}[ht]		
\centering		
    \begin{subfigure}{0.45\textwidth}
		\centering			
        \def\svgwidth{4.4cm}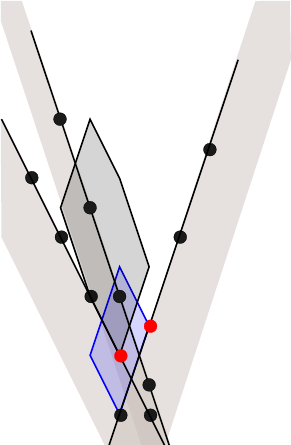
		\caption{The parallelogram $\mathcal{Q}$ (shaded) used to locate the translate $\mathcal{K}$ of $\mathcal{S}$. The point $u'$ marks the corner of $\st_3(U,V)$ that falls outside $\st_2(U,V)$, and $u$ is the anchor point on the $(-\pmb{\ell}_2)$-boundary of $\st_2(U,V)$ where the edge $\mathcal{K} \cap (\pmb{\ell}_2)_{\mathcal{K}}$} terminates.
		\end{subfigure}
		\hfill
		\begin{subfigure}{0.45\textwidth}
		\centering
		\def\svgwidth{4.4cm}
\begingroup%
  \makeatletter%
  \providecommand\color[2][]{%
    \errmessage{(Inkscape) Color is used for the text in Inkscape, but the package 'color.sty' is not loaded}%
    \renewcommand\color[2][]{}%
  }%
  \providecommand\transparent[1]{%
    \errmessage{(Inkscape) Transparency is used (non-zero) for the text in Inkscape, but the package 'transparent.sty' is not loaded}%
    \renewcommand\transparent[1]{}%
  }%
  \providecommand\rotatebox[2]{#2}%
  \newcommand*\fsize{\dimexpr\f@size pt\relax}%
  \newcommand*\lineheight[1]{\fontsize{\fsize}{#1\fsize}\selectfont}%
  \ifx\svgwidth\undefined%
    \setlength{\unitlength}{139.71926592bp}%
    \ifx\svgscale\undefined%
      \relax%
    \else%
      \setlength{\unitlength}{\unitlength * \real{\svgscale}}%
    \fi%
  \else%
    \setlength{\unitlength}{\svgwidth}%
  \fi%
  \global\let\svgwidth\undefined%
  \global\let\svgscale\undefined%
  \makeatother%
  \begin{picture}(1,1.52781305)%
    \lineheight{1}%
    \setlength\tabcolsep{0pt}%
    \put(0,0){\includegraphics[width=\unitlength,page=1]{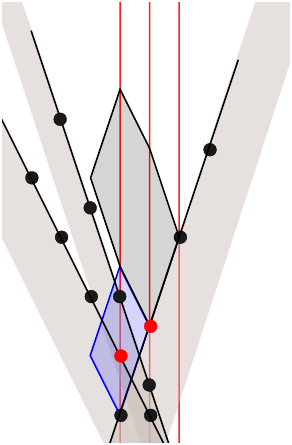}}%
    \put(0.86373336,1.11425185){\color[rgb]{0,0,0}\transparent{0.37918204}\rotatebox{74.39525}{\makebox(0,0)[lt]{\lineheight{1.25}\smash{\begin{tabular}[t]{l}$\st_3(U,V)$\end{tabular}}}}}%
    \put(0.01765124,1.47317305){\color[rgb]{0,0,0}\transparent{0.37918204}\rotatebox{-70.910808}{\makebox(0,0)[lt]{\lineheight{1.25}\smash{\begin{tabular}[t]{l}$\st_1(U,V)$\end{tabular}}}}}%
    \put(0.0257807,0.83400572){\color[rgb]{0,0,0}\transparent{0.37918204}\rotatebox{-65.344665}{\makebox(0,0)[lt]{\lineheight{1.25}\smash{\begin{tabular}[t]{l}$\st_2(U,V)$\end{tabular}}}}}%
    \put(0.24907004,1.26183279){\color[rgb]{0,0,0}\makebox(0,0)[lt]{\lineheight{1.25}\smash{\begin{tabular}[t]{l}$\mathcal{K}+u'-u$\end{tabular}}}}%
  \end{picture}%
\endgroup%
			
        \caption{The translate $\mathcal{K}+u'-u$ inside $U_1$. The three vertical lines (red), counted from the initial vertex of $-w_1$ to the final vertex of $w_1$ inclusive, reflect the width constraint imposed by $\mathcal S$ being contained in the rectangle $R_{4,n}$.}
		\end{subfigure}
		
	\caption{Case~1(a): construction of the translate $\mathcal{K}$ of
$\mathcal{S}$ used to establish $\pmb{\ell}_2$-compatibility of $W''$.}
	\label{fig5}	
\end{figure}

Thus, since $\mathcal{K}+u'-u \subset U_1$, every
non-vertical edge has exactly two integer points, and each direction
vector $v_{\pmb{\ell}_i}$ advances from one vertical line to the
adjacent one, it follows that
\[
\sst_{\mathcal{K}+u'-u-v_{\pmb{\ell}_3}
-v_{\pmb{\ell}_2}}(-\pmb{\ell}_2) \cap \st_1(U,V) \neq \emptyset.
\]
A fortiori, since $u \in \mathcal{Q}$, we also have
$\sst_{\mathcal{K}}(-\pmb{\ell}_2) \cap \st_1(U,V) \neq \emptyset$.
Therefore, $W''$ is $\pmb{\ell}_2$-compatible with some translate of $\mathcal{S}$, 
and Proposition~\ref{prop_CyrKra}(i) contradicts Theorem~\ref{thm:main general theorem}.

\medskip

\textit{Case 1(b): the $y$-axis is nonexpansive on $X_{\vartheta}$.}

\medskip

We may assume $\pmb{\ell}_1$ is parallel to the oriented $y$-axis. The
orientations of $\pmb{\ell}_2$ and $\pmb{\ell}_3$, and hence the
strip-free regions $W$, $W'$ and $W''$, are kept exactly as fixed in
Case~1(a) (see~\eqref{eq:strip-free-triple-m3} above). Among the edges of
$\mathcal{S}$, we single out the six that come in antiparallel pairs
directed along the three nonexpansive lines: denote them
$-\varpi_1,\varpi_1,-\varpi_2,\varpi_2,-\varpi_3,\varpi_3$, with
$\varpi_i$ parallel to $\pmb{\ell}_i$, and with $\varpi_1 \prec \varpi_2
\prec \varpi_3$ in the successor order of $E(\mathcal{S})$. Any remaining
edges of $\mathcal{S}$ are directed along expansive lines and play no
role in the argument. Two edge-cardinality observations will drive the
argument:
\begin{enumerate}[label=(\roman*)]\setlength{\itemsep}{5pt}
\item Since $\mathcal{S} \subset R_{4,n}$, at most one non-vertical
edge $w \in E(\mathcal{S})$ between $-\varpi_1$ and $\varpi_1$ can
satisfy $|w \cap \mathcal{S}| = 3$.

\item If $\varpi_i$ ($i \neq 1$) satisfies $|\varpi_i \cap \mathcal{S}| = 3$,
then every oriented line $\pmb{\ell}'_i$ parallel to $\varpi_i$ with
$\pmb{\ell}'_i \cap \mathcal{S} \neq \emptyset$ satisfies
$|\pmb{\ell}'_i \cap \mathcal{S}| \geq 2$. Indeed, since $\mathcal{S}$
is convex, the integer points of $-\varpi_i \cap \mathcal{S}$ can be
translated vertically within $\mathcal{S}$ and then, if necessary,
shifted by $\pm v_{\pmb{\ell}_i}$ within $\mathcal{S}$ to reach
$\varpi_i$. (This ensures that, for $\pmb{\ell} = \pmb{\ell}_i$ in Proposition~\ref{prop_CyrKra}, 
we have \( B = \mathcal{S} \setminus (\pmb{\ell}_i)_{\mathcal{S}} \).)
\end{enumerate}

In what follows we do not distinguish between edges of $\mathcal{S}$ and
those of a translate of $\mathcal{S}$.

\smallskip
We first ask whether $W'$ is $\pmb{\ell}_3$-compatible with some translate
$\mathcal{T}$ of $\mathcal{S}$. Suppose that such a translate exists;
by shifting $\mathcal{T}$ along $\pmb{\ell}_3$ if necessary, we may
position it so that $W'$ is $\pmb{\ell}_3$-compatible with $\mathcal{T}$
but not with $\mathcal{T}+v_{\pmb{\ell}_3}$, that is, $\mathcal{T}$
is as far as possible in the direction of $\pmb{\ell}_3$ while still
satisfying the compatibility condition. There are then two possibilities to consider:

\smallskip
\textit{Subcase $|\varpi_3 \cap \mathcal{T}| = 2$.} Then
Proposition~\ref{prop_CyrKra}(i) gives $\ell_3$-periodicity of $\vartheta$
on $\sst_{\mathcal{T}\setminus\varpi_3}(\pmb{\ell}_3)$, contradicting
Theorem~\ref{thm:main general theorem}.

\smallskip
\textit{Subcase $|\varpi_3 \cap \mathcal{T}| = 3$ (so $|\varpi_2 \cap
\mathcal{T}| = 2$).} The assumption on
$\mathcal{T}+v_{\pmb{\ell}_3}$ forces
\begin{equation}\label{eq:position of T - caso m=3(b)}
    (\pmb{\ell}_1)_{\mathcal{T}} \cap \bb{Z}^2
    = (-\pmb{\ell}_{1})_{\st_1(U,V)} \cap \bb{Z}^2.
\end{equation}
Is $W'$ also $\pmb{\ell}_3$-compatible with $\mathcal{T}-2v_{\pmb
{\ell}_3}$? If so, Proposition~\ref{prop_CyrKra}(ii) yields $\ell_3$-periodicity,
again contradicting Theorem~\ref{thm:main general theorem}.

Otherwise, let $ g $ be the initial point of the edge $\varpi_3 \in E(\mathcal{T})$ and define
the parallelogram
\[
\mathcal{Q} := \mathrm{Conv}\!\left\{
g,\;
g-2v_{\pmb{\ell}_3},\;
g-v_{\pmb{\ell}_2},\;
g-2v_{\pmb{\ell}_3}-v_{\pmb{\ell}_2}
\right\} \cap \bb{Z}^2.
\]
Because $\mathcal{T}$ straddles four vertical lines between $-\varpi_1$
and $\varpi_1$ (inclusive), equation~\eqref{eq:position of T - caso m=3(b)} gives
\[
\sst_{\mathcal{T}-2v_{\pmb{\ell}_3}-v_{\pmb{\ell}_2}}(-\pmb{\ell}_2)
\cap \st_1(U,V) \neq \emptyset.
\]
Choosing $g' \in \mathcal{Q} \cap (-\pmb{\ell}_2)_{\st_2(U,V)}$
with $g' \notin \st_3(U,V)$ and setting $\mathcal{U} = \mathcal{T}
+g'-g$, the region $W''$ is $\pmb{\ell}_2$-compatible with
$\mathcal{U}$ (see Figure~\ref{fig6}). By Proposition~\ref{prop_CyrKra}(i), $\vartheta\sob{\sst_{\mathcal{U} \setminus\varpi_2}(-\pmb{\ell}_2)}$ is $\ell_2$-periodic, contradicting Theorem~\ref{thm:main general theorem}.

\begin{figure}[ht]
    \centering
    \def\svgwidth{5cm}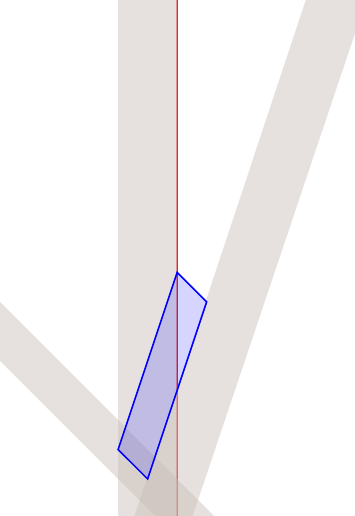
    \caption{The parallelogram $\mathcal{Q}$ (shaded) and the translate $\mathcal{U}$ of $\mathcal{S}$ positioned so that $W''$ is
    $\pmb{\ell}_2$-compatible with $\mathcal{U}$. The point $g$ is the initial point of $\varpi_3 \in E(\mathcal{T})$, and $g'$ is chosen on the $(-\pmb{\ell}_2)$-boundary of $\st_2(U,V)$ just outside $\st_3(U,V)$.}
    \label{fig6}
\end{figure}

\smallskip
It remains to handle the case where $W'$ is not $\pmb{\ell}_3$-compatible
with \emph{any} translate of $\mathcal{S}$. The argument now mirrors
Case~1(a): define $u'$, $\mathcal{Q}$, $u$, and $\mathcal{K}$
exactly as there (see Figure \ref{fig7}{\color{red}(A)}). As before the translate $\mathcal{K}+u'-
u$ lies in $U_1$, and now it forces five vertical lines between
$(-\pmb{\ell}_1)_{\mathcal{K}+u'-u}$ and
$(-\pmb{\ell}_1)_{\st_1(U,V)}$, inclusive
(see Figure~\ref{fig7}{\color{red}(B)}). This gives
\[
\sst_{\mathcal{K}+u'-u+v_{\pmb{\ell}_3}-3v_{\pmb{\ell}_2}}(-\pmb{\ell}_2) \cap \st_1(U,V) \neq  \emptyset,
\] 
or equivalently,
\[(
\sst_{\mathcal{K}+u'-u+v_{\pmb{\ell}_3}-v_{\pmb{\ell}_2}}(-\pmb{\ell}_2)-2v_{\pmb{\ell}_2}) \cap \st_1(U,V) \neq  \emptyset.
\]
Since $u \in \mathcal{Q}$, in both cases
($|\varpi_2 \cap \mathcal{K}|=2$ or $|\varpi_2 \cap \mathcal{K}|=3$)
Proposition~\ref{prop_CyrKra} supplies $\ell_2$-periodicity on the
appropriate semi-strip, contradicting Theorem~\ref{thm:main general theorem}.

\begin{figure}[ht]
	\begin{minipage}[t]{0.45\linewidth}
		\begin{center}
			\def\svgwidth{4.4cm}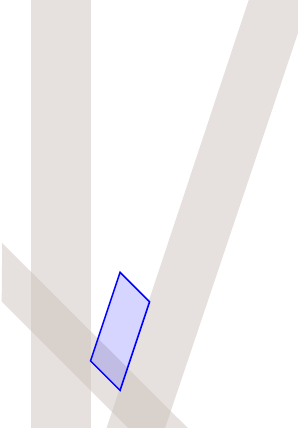\\
		\end{center}
		\small{(A) The parallelogram $\mathcal{Q}$ (shaded) and the anchor
points $u'$ and $u$ used to locate the translate $\mathcal{K}$ of
$\mathcal{S}$, defined exactly as in Case~1(a). The point $u'$ marks
the corner of $\st_3(U,V)$ falling outside $\st_2(U,V)$, and $u$ is
chosen on the $(-\pmb{\ell}_2)$-boundary of $\st_2(U,V)$ just outside
$\st_3(U,V)$, so that the edge $\mathcal{K} \cap (\pmb{\ell}_2)_{\mathcal{K}}$ terminates at~$u$.}
	\end{minipage} 
	\hfill
	\begin{minipage}[t]{0.45\linewidth}
		\begin{center}
			\def\svgwidth{4.4cm}
\begingroup%
  \makeatletter%
  \providecommand\color[2][]{%
    \errmessage{(Inkscape) Color is used for the text in Inkscape, but the package 'color.sty' is not loaded}%
    \renewcommand\color[2][]{}%
  }%
  \providecommand\transparent[1]{%
    \errmessage{(Inkscape) Transparency is used (non-zero) for the text in Inkscape, but the package 'transparent.sty' is not loaded}%
    \renewcommand\transparent[1]{}%
  }%
  \providecommand\rotatebox[2]{#2}%
  \newcommand*\fsize{\dimexpr\f@size pt\relax}%
  \newcommand*\lineheight[1]{\fontsize{\fsize}{#1\fsize}\selectfont}%
  \ifx\svgwidth\undefined%
    \setlength{\unitlength}{142.66671609bp}%
    \ifx\svgscale\undefined%
      \relax%
    \else%
      \setlength{\unitlength}{\unitlength * \real{\svgscale}}%
    \fi%
  \else%
    \setlength{\unitlength}{\svgwidth}%
  \fi%
  \global\let\svgwidth\undefined%
  \global\let\svgscale\undefined%
  \makeatother%
  \begin{picture}(1,1.43691505)%
    \lineheight{1}%
    \setlength\tabcolsep{0pt}%
    \put(0,0){\includegraphics[width=\unitlength,page=1]{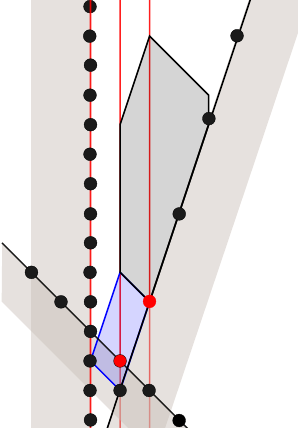}}%
    \put(0.34527326,1.36140538){\color[rgb]{0,0,0}\makebox(0,0)[lt]{\lineheight{1.25}\smash{\begin{tabular}[t]{l}$\mathcal{K}+u'-u$\end{tabular}}}}%
    \put(0.0065499,0.48089584){\color[rgb]{0.16078431,0.14901961,0.14117647}\transparent{0.51354867}\rotatebox{-44.618246}{\makebox(0,0)[lt]{\lineheight{1.25}\smash{\begin{tabular}[t]{l}$\st_2(U,V)$\end{tabular}}}}}%
    \put(0.81576276,1.03709617){\color[rgb]{0.16078431,0.14901961,0.14117647}\transparent{0.51378149}\rotatebox{70.72953}{\makebox(0,0)[lt]{\lineheight{1.25}\smash{\begin{tabular}[t]{l}$\st_3(U,V)$\end{tabular}}}}}%
    \put(0.18635003,1.39795882){\color[rgb]{0.16078431,0.14901961,0.14117647}\transparent{0.51378149}\rotatebox{-90}{\makebox(0,0)[lt]{\lineheight{1.25}\smash{\begin{tabular}[t]{l}$\st_1(U,V)$\end{tabular}}}}}%
    \put(0,0){\includegraphics[width=\unitlength,page=2]{fig7b.pdf}}%
  \end{picture}%
\endgroup%
\\
		\end{center}
		\small{(B) The translate $\mathcal{K}+u'-u$ inside $U_1$. The five vertical lines (red) between $(-\pmb{\ell}_1)_{\mathcal{K}+u'-u}$ and the $(-\pmb{\ell}_1)$-boundary of $\st_1(U,V)$ reflect the additional column available compared to Case~1(a), and are the key to reaching $\st_1(U,V)$ via the $\pmb{\ell}_3$- and $\pmb{\ell}_2$-shifts.}
	\end{minipage}
	\caption{Case~1(b), subcase where $W'$ is not $\pmb{\ell}_3$-compatible with any translate of $\mathcal{S}$. The construction mirrors Case~1(a), but the translation analysis of $\mathcal{K}+u'-u$ now takes into account five vertical lines (shown in red in ~(B)), which is precisely what allows the combined $\pmb{\ell}_3$- and $\pmb{\ell}_2$-shifts to reach $\st_1(U,V)$ and confirm $\pmb{\ell}_2$-compatibility of $W'''$.}
	\label{fig7}
\end{figure}

\medbreak

\textbf{Case 2: $m = 4$.}

\medbreak

When $m = 4$, the constraint $\mathcal{S} \subset R_{4,n}$ is more
stringent: any extremal $\vartheta$-generating set in $R_{4,n}$ has exactly four
pairs of antiparallel edges, two of which are necessarily directed along the $y$-axis, 
and \emph{every} non-vertical edge is forced to contain exactly two integer points. 
After renaming if necessary, we orient $\ell_1$ along the $y$-axis and index $\ell_2,\ell_3,\ell_4$ counterclockwise from
$\pmb{\ell}_1$, choosing their orientations so that the following four strip-free regions are well defined: 
\begin{equation*}
\left.
\begin{aligned}
&W &&\text{an } (-\pmb{\ell}_4,\pmb{\ell}_3)\text{-strip-free region for }(U,V),\\
&W' &&\text{an } (\pmb{\ell}_1,\pmb{\ell}_4)\text{-strip-free region for }(U,V),\\
&W'' &&\text{an } (\pmb{\ell}_2,-\pmb{\ell}_1)\text{-strip-free region for }(U,V), \\
&W''' &&\text{an } (\pmb{\ell}_3,-\pmb{\ell}_2)\text{-strip-free region for }(U,V).
\end{aligned}
\right.
\end{equation*}

Is $W'$ $\pmb{\ell}_4$-compatible with some translate $\mathcal{T}$ of $\mathcal{S}$? If so, since every non-vertical edge of $\mathcal{T}$ contains exactly two integer points, Proposition~\ref{prop_CyrKra}(i) immediately yields $\ell_4$-periodicity of $\vartheta$ on $\sst_{\mathcal{T}\setminus (\pmb{\ell}_4)_{\mathcal{T}}}(\pmb{\ell}_4)$, contradicting Theorem~\ref{thm:main general theorem}.
Hence, $W'$ is not $\pmb{\ell}_4$-compatible with any translate of $\mathcal{S}$. Is $W''$ $\pmb{\ell}_2$-compatible with some translate of $\mathcal{S}$? If so, then, as before, we obtain a contradiction with Theorem~\ref{thm:main general theorem}. Hence, $W''$ is not $\pmb{\ell}_2$-compatible with any translate of $\mathcal{S}$.
It remains to examine what the $\pmb{\ell}_2$-incompatibility of $W''$ concretely entails in this situation. 
Since the $\pmb{\ell}_4$-incompatible has already been verified and every non-vertical edge of $\mathcal{S}$ contains exactly two integer points, the only geometric obstruction is that for any translate $\mathcal{T}$ of $\mathcal{S}$ such that 
$\mathcal{T} \cap (W'')^1 \neq \emptyset $, the condition $\sst_{\mathcal{T}}(-\pmb{\ell}_2) \cap \st_1(U,V) \neq \emptyset$ means that $\mathcal{T}$ meets $\st_3(U,V)$. 
In other words, part of $\st_3(U,V)$ is positioned far enough from $\st_1(U,V) \cap \st_2(U,V) $ inside $ \mathcal H(\pmb{\ell}_2) \cap \mathcal H(-\pmb{\ell}_1) $ that no such translate of $\mathcal{S}$ can simultaneously reach $\st_2(U,V)$ and have its $(-\pmb{\ell}_2)$-semi-strip enter $\st_1(U,V)$ without invading $\st_3(U,V)$.
This separation of $\st_3(U,V)$ with respect to $\st_1(U,V) \cap \st_2(U,V) $ is precisely what allows
$W'''$ to be $\pmb{\ell}_3$-compatible: there exists a translate $\mathcal{T}$ of $\mathcal{S}$ meeting $\st_3(U,V)$ with its remainder $\mathcal{T}\setminus(\pmb{\ell}_3)_{\mathcal{T}}$ contained in $(W''')^2$,
while its $(-\pmb{\ell}_3)$-semi-strip enters $\st_2(U,V)$, so that both
conditions of Definition~\ref{def:compatible} are satisfied.
But then Proposition~\ref{prop_CyrKra}(i) yields $\ell_3$-periodicity on a semi-strip, again contradicting Theorem~\ref{thm:main general theorem}.

In every case we reach a contradiction, so $\eta$ must be periodic.
\end{proof}

We highlight that the case in which \(P_{\eta}(5,n) \leq 5n\) for some \(n \in \mathbb{N}\) cannot be handled using the same method. In this case, the geometry of generating sets contained in \(R_{5,n} = \llb 5 \rrb \times \llb n \rrb\) is richer, allowing a generating set with a non-vertical pair of antiparallel edges with discrepancy two. Consequently, for such a generating set, the set \(B\) in Proposition~\ref{prop_CyrKra} may not be large enough to yields a contradiction in Theorem~\ref{thm:main general theorem}.

\end{document}